\documentclass[12pt]{article}
\usepackage{amsthm,amsmath,amssymb,amscd,amsfonts,color,graphicx,bm}
\newcommand{\intprod}{\mathbin{\raisebox{\depth}{\scalebox{1}[-1]{$\lnot$}}}}
\begin{document}

\centerline{\bf New Procedure to Generate Multipliers in Complex}
\centerline{\bf Neumann Problem and Effective Kohn Algorithm}
\bigbreak\centerline{\it Dedicated to the Memory of Professor Lu Qikeng}

\bigbreak
\centerline{Yum-Tong Siu}

\bigbreak

The purpose of this note is threefold.  (i) To explain the effective Kohn algorithm for multipliers in the complex Neumann problem and its difference with the full-real-radical Kohn algorithm, especially in the context of an example
of Catlin-D'Angelo concerning the ineffectivness of the latter.  (ii) To extend the techniques of
multiplier ideal sheaves for the complex Neumann problem to general systems of partial differential equations. (iii) To present a new procedure of generation of
multipliers in the complex Neumann problem as a special case of the multiplier ideal sheaves techniques for general systems of partial differential equation.

\medbreak For {\it a priori} estimates in the theory of partial differential equations some of the standard techniques are the following. (i) Using integration by parts to get $L^2$ estimates of derivatives in certain directions, for example, $L^2$ estimates of all first-order partial derivatives of a function with compact support from applying integration by parts to its inner product with its Laplacian. (ii) Using Lie bracket of two vector fields to conclude, from given derivative estimates of fractional orders along each of the two vectors, the derivative estimates of lower fractional order along their Lie bracket, for example, H\"ormander's work on the sum of squares of vector fields [H\"ormander1967].

\medbreak The technique of multiplier ideal sheaves introduced by Kohn for the complex Neumann problem is a new method to conclude, from $L^2$ estimates of derivatives along certain {\it complex-valued} vector fields, the derivative estimates of lower fractional orders in all directions by introducing the notion of multipliers.  In a system of partial differential equations where the estimate is for a vector-valued test function with components $\psi_\nu$, for given complex-valued vector fields $Y_j$, when the estimates for several linear combinations $\sum_{j,\nu}\rho_{j,\nu}Y_j\psi_\nu$ with some given smooth functions $\rho_{j,\nu}$ are known, the multipliers are
the smooth coefficients $a_\nu$ of a linear combination $\sum_\nu a_\nu\psi_\nu$ such that there is an estimate of the Sobolev $L^2$ norm of $\sum_\nu a_\nu\varphi_\nu$ of some positive fractional order.  The goal is to derive differential relations among the multipliers to obtain some geometric condition to solve the regularity problem for the given system of partial differential equations, by the method of using the differential relations and some initial multipliers to conclude that the function which is identically $1$ is a scalar multiplier.

\medbreak For Kohn's original setting of a bounded weakly pseudoconvex domain $\Omega$ in ${\mathbb C}^n$ the vector fields $Y_j$ are the vector fields $L_1,\cdots,L_{n-1}$ of type $(1,0)$ tangential to the boundary $\partial\Omega$ of $\Omega$ together with their complex conjugates $\bar L_1,\cdots,\bar L_{n-1}$.  The test function is a $(0,1)$-form in the domains of $\bar\partial$ and $\bar\partial^*$ whose $n-1$ tangential components are $\varphi_{\overline 1},\cdots,\varphi_{\overline{n-1}}$.  The given linear combinations $\sum_{j,\nu}\rho_{j,\nu}Y_j\psi_\nu$ are the $(n-1)^2+1$ linear combinations $\bar L_j\varphi_\nu$ (for $1\leq j,\nu\leq n-1$) and $\sum_{\nu=1}^{n-1}L_\nu\varphi_\nu$.

\medbreak The structure of this note is as follows.  We start out with the background and motivation for the technique of multiplier ideal sheaves.  We then discuss the two Kohn algorithm of generating multiplies.  One is the full-real-radical Kohn algorithm.  The other is the effective Kohn algorithm.  We explain the algebraic geometric techniques in the effective Kohn algorithm, with special attention paid to the effectiveness of the orders of subellipticity in each step.  The effective Kohn algorithm is then applied to Catlin-D'Angelo's example to highlight the difference in effectiveness between the full-real-radical Kohn algorithm and the effective Kohn algorithm.  We present the generalization of the technique of multipliers to general systems of partial differential equations.  Finally we apply the generalized techniques to the complex Neumann problem to obtain a new procedure to generate vector multipliers from matrix multipliers for special domains.

\medbreak The notations ${\mathbb N}$, ${\mathbb R}$, and ${\mathbb C}$ mean respectively all positive integers, all real numbers, and all complex numbers.  The notations ${\mathcal O}_{{\mathbb C}^n,P}$ and ${\mathfrak m}_{{\mathbb C}^n,P}$ mean respectively all holomorphic germs on ${\mathbb C}^n$ at $P$ and the maximum ideal of the local ring ${\mathcal O}_{{\mathbb C}^n,P}$.  Unless specified otherwise, $\left\|\cdot\right\|$ denotes the $L^2$ norm and $\left(\cdot,\,\cdot\right)$ denotes the $L^2$ inner product.  The notation $\left\|\cdot\right\|_{L^2(U)}$ is also used to more clearly specify that it is the $L^2$ norm over $U$.  The notation $\left(\cdot,\,\cdot\right)_{L^2(U)}$ is also used to more clearly specify that it is the $L^2$ inner product over $U$.  The notation $L^2_k$ means the Sobolev norm defined by using the $L^2$ norm of derivatives up to order $k$.  The notations $\partial_j$ and $\bar\partial_j$ mean respectively $\frac{\partial}{\partial z_j}$ and $\frac{\partial}{\partial\bar z_j}$, where $z_1,\cdots,z_n$ are the coordinates of ${\mathbb C}^n$.

\vskip.3in\noindent{\sc Table of Contents}

\medbreak\noindent\S 1. Background and Motivation for Multiplier Ideal Sheaves

\medbreak\noindent\S 2. Generation of Multipliers in Kohn's Algorithm

\medbreak\noindent\S 3. Orders of Subellpiticity in Algebraic Geometric Techniques for 3-Dimensional Special Domain

\medbreak\noindent\S 4.  Effective Kohn Algorithm Applied to Catlin-D'Angelo's Example

\medbreak\noindent\S 5. Multipliers in More General Setting

\medbreak\noindent\S 6. New Procedure to Generate Vector Multiplier from Matrix Multiplier in Complex Neumann Problem of Special Domain

\vskip.5in\noindent{\bf\S 1.} {\sc Background and Motivation for Multiplier Ideal Sheaves}

\bigbreak In regularity problems of local or global systems of partial differential equations, multiplier ideal sheaves describe the location and the extent of the failure of {\it a priori} estimates.  There are two ways to introduce such multiplier ideal sheaves.  The first way is from the continuity method of solving partial differential equations (which are usually nonlinear partial differential equations defined on compact Riemannian manifolds) where a multiplier ideal sheaf arises as the limit of rescaling of local coordinate charts to make possible the use of Ascoli-Arzela techniques for convergence.  The second way is just to directly introduce multiplier ideal sheaves as factors required in the integral norms to make {\it a priori} estimates hold for the regularity problem.  We very briefly describe both.

\bigbreak\noindent(1.1) {\it Multiplier Ideal Sheaves from Limit of Rescaling of Local Coordinate Charts.}
One method to solve nonlinear partial differential equations is to use the {\it continuity method} which uses a family of partial differential equations
parametrized by $t\in[0,1]$ so that
\begin{itemize}
\item[(i)] $t=1$ is the given partial differential equations,
\item[(ii)] $t=0$ can be solved (which means solvability of partial differential equation at initial parameter), and
\item[(iii)] from the solution for $t=t_0<1$ it is possible to solve for $t<t_0+\varepsilon$ for some $\varepsilon>0$ (which means that the {\it openness} property
is assumed).
\end{itemize}
The difficulty is to prove the {\it closedness} property that for $0<t^*\leq 1$ solutions $s_\nu$ for the parameter value $t=t_\nu$ for $t_\nu\nearrow t^*$ can be used to construct a solution $s$ for the parameter value $t=t^*$.

\medbreak The natural approach is to obtain $s$ by taking the limit of $s_{\nu_k}$ for some subsequence $\{\nu_k\}$ of the sequence $\{\nu\}$.  Usually the boundedness in some weak norm for $s_\nu$, for example,
$L^2$, can be derived from the setup of the given partial differential equation.  The method of Ascoli-Arzela calls for boundedness in some stronger norm for $s_\nu$, for example, $L_1^2$ defined
by the first-order derivatives being $L^2$.

\medbreak
In the setting of a manifold the derivative involved in $L^2_1$ depends on the choice of local coordinates so that $L^2_1$ is determined up to equivalence (of sandwiching
between its products with two constants).  One can always cheat by using different rescaling of local coordinates for each $s_\nu$, but one has to pay the price that at the end
the variable rescaling of local coordinates results in a factor (or {\it multiplier}) which is the limit of the Jacobian determinant in coordinate change in the integral for $L^2$.

\smallbreak The limit $s$ of $s_{\nu_k}$ from the Ascoli-Arzela argument is $L^2$ only after the insertion of the multiplier in its $L^2$ integral.  This can be interpreted as transforming the requirement in the the method of Ascoli-Arzela for two norms $L^2$ and $L_1^2$ (with difference in orders of differentiation in their definitions) to two norms which are the $L^2$ norm and the multiplier-modified $L_1^2$ norm.  Multipliers form an ideal sheaf, called the {\it multiplier ideal sheaf}.  Global conditions, such as topological conditions, can be used to conclude that the multiplier ideal sheaf must be the full structure sheaf, giving the solvability of the partial differential equation from global conditions.  Examples are (i) the existence of Hermitian-Einstein metrics for stable holomorphic vector bundles over compact algebraic (or K\"ahler) manifolds, where the method of limit of rescaling of local coordinates was first applied by Donaldson in [Donaldson1985] for the surface case and (ii) the existence of K\"ahler-Einstein metrics for certain Fano manifolds, where the method of multiplier ideal sheaves as defined by taking limit of metrics of the anti-canonical line bundle was first applied by Nadel in [Nadel1990].

\bigbreak We now look at the second way of introducing multiplier ideal sheaves, which is best described by using, as an example, the problem of the regularity of the Kohn solution for the complex Neumann problem on a bounded weakly pseudoconvex domain in ${\mathbb C}^n$ with smooth boundary.

\bigbreak\noindent(1.2) {\it Regularity Problem of Kohn Solution for Complex Neumann Problem on Weakly Pseudoconvex Domain.}  Let $\Omega$ be a bounded weakly pseudoconvex domain in ${\mathbb C}^n$ with smooth boundary $\partial\Omega$.  It means that there is a smooth function $r$ defined on some open neighborhood $W$ of $\partial\Omega$ in ${\mathbb C}^n$ such that

\medbreak\noindent(i)  $\Omega\cap W=\left\{\,z\in W\,\big|\,r<0\,\right\}$,

\medbreak\noindent(ii) $dr$ is nowhere zero on $\partial\Omega$, and

\medbreak\noindent(iii) the $(1,1)$-form $\sqrt{-1}\,\partial\bar\partial r$ on $W$ assumes nonnegative value when evaluated at $(\xi,\bar\xi)$ for $\xi\in T^{(1,0)}_{\partial\Omega}$, where $T^{(1,0)}_{\partial\Omega}$ is the bundle of all tangent vectors $\xi$ of ${\mathbb C}^n$ of type $(1,0)$ at points of $\partial\Omega$ which are tangential to $\partial\Omega$ (in the sense that $\xi(r)=0$).

\bigbreak For notational simplicity we will consider only the complex Neumann problem in the case of $(0,1)$-forms (instead of $(0,p)$-form for $1\leq p\leq n$).  For a $\bar\partial$-closed $(0,1)$-form $f$ on $\Omega$ which are smooth on the closure $\bar\Omega$ of $\Omega$, the {\it Kohn solution} $u$ for the $\bar\partial$-equation $\bar\partial u=f$ on $\Omega$ means the unique smooth function $u$ on $\Omega$ such that $\bar\partial u=f$ on $\Omega$ and $u$ is perpendicular to all $L^2$ holomorphic functions on $\Omega$ with respect to the usual Euclidean volume form of ${\mathbb C}^n$.

\medbreak The regularity problem of the Kohn solution for the complex Neumann problem on a bounded weakly pseudoconvex domain with smooth boundary is to study under what additional assumption on $\Omega$ the Kohn solution $u$ is always smooth on $\bar\Omega$ when the given $(0,1)$-form $f$ is smooth on $\bar\Omega$.

\bigbreak\noindent(1.3) {\it Regularity from Subelliptic Estimate.}  The {\it subelliptic estimate} of order $\varepsilon>0$ is said to hold at a point $P$ of $\partial\Omega$ if there exist an
open neighborhood $U$ of $P$ in ${\mathbb C}^n$ and positive numbers
$\varepsilon$ and $C$ satisfying
$$
\||\varphi|\|_\varepsilon^2\leq C\left(\|\bar\partial\varphi\|^2+\|\bar\partial^*
\varphi\|^2+\|\varphi\|^2\right)
$$
for all smooth $(0,1)$-forms $\varphi$ on $U\cap\bar\Omega$ with compact support which belong to
the domain of the {\it actual} adjoint $\bar\partial^*$ of $\bar\partial$ (with respect to the usual $L^2$ inner product), where

\medbreak\noindent(i) $\||\cdot|\|_\varepsilon$ is the Sobolev $L^2$ norm on $\Omega$
involving derivatives up to order $\varepsilon$ in the boundary
tangential directions of $\Omega$, and
\medbreak\noindent(ii) $\|\cdot\|$ is the usual $L^2$ norm on $\Omega$ without involving any
derivatives.

\bigbreak See [Kohn1979, p.92, (3.4)] for a detailed definition of the Sobolev norm $\||\cdot|\|_\varepsilon$.  In this note we will also use the notation $\Lambda^s$ introduced in [Kohn1979, p.92, (3.3)] which is the pseudo-differential operator corresponding to the $\left(\frac{s}{2}\right)$-th power of $1$ plus the Laplacian in tangential coordinates of the boundary of $\Omega$.

\bigbreak We would like to comment on the reason for the use of the Sobolev norm $\||\cdot|\|_\varepsilon$ in whose definition only
derivatives along the tangent directions of $\partial\Omega$ is used.  The condition for a smooth $(0,1)$-form $g$ on $\bar\Omega$ to belong to the domain of the {\it actual} adjoint
$\bar\partial^*$ of $\bar\partial$ (instead of the {\it formal} adjoint of $\bar\partial$) is that the components of $g$ normal to $\partial\Omega$ vanish at all points of $\partial\Omega$ (which,
in other words, means that the pointwise inner product of $g$ and $\bar\partial r$ vanishes at every point of $\partial\Omega$).  When $g$ is in the domain of $\bar\partial^*$, the derivative of $g$ along tangent directions of $\partial\Omega$ still belongs to the domain of $\bar\partial^*$, but if $g$ is differentiated in the normal direction of $\partial\Omega$, the result in general will no longer belong to the domain of $\bar\partial^*$.  That is the reason why we would like to avoid using differentiation along the normal directions of $\bar\partial$ in the Sobolev norm of order $\varepsilon$ adopted in the definition of the subelliptic estimate.

\medbreak Sobolev norms involving derivatives are used to enable us to conclude the order of differentiability of the solution $u$ of the $\bar\partial$-equation $\bar\partial u=f$ from the order of differentiability of the right-hand side $f$ of the equation.  Though we do not include the differentiation in the normal direction of $\partial\Omega$ in the Sobolev norm used, from the order of differentiability of the solution $u$ along the {\it tangent directions} of $\partial\Omega$ we can still conclude the order of differentiability of $u$ along the normal direction of $\partial\Omega$ because the equation $\bar\partial u=f$ itself provides us directly the differentiability of $u$ along the {\it real normal} direction of $\partial\Omega$ from the differentiability of $u$ along the $(0,1)$ component of the {\it complex normal} direction of $\partial\Omega$.

\medbreak It was proved by Kohn-Nirenberg in 1965 [Kohn-Nirenberg1965, p.458, Theorem 4] that if for some $\varepsilon>0$ the subelliptic estimate of order $\varepsilon$ holds at every point of $\partial\Omega$, then the smoothness on $\bar\Omega$ of the Kohn solution $u$ of $\bar\partial u=f$ follows from the smoothness of the $\bar\partial$-closed $(0,1)$-form $f$ on $\bar\Omega$.

\bigbreak\noindent(1.4) {\it Multipliers to Measure Location and Extent of Failure of A Priori Estimates}.  A smooth function germ $F$ at a point $P$ of $\partial\Omega$ (defined on some open neighborhood $U_F$ of $P$ in ${\mathbb C}^n$ is called a {\it scalar multiplier} if for some positive number $\varepsilon_F>0$ and some positive constant $C_{{}_F}$ the subelliptic estimate
$$
\||F\varphi|\|_{\varepsilon_{{}_F}}^2\leq C_{{}_F}\left(\|\bar\partial
\varphi\|^2+\|\bar\partial^*\varphi\|^2+\|\varphi\|^2\right)\leqno{(1.4.1)}
$$
of order $\varepsilon_F$, {\it modified} by the factor $F$, holds for
for all smooth $(0,1)$-forms $\varphi$ on $U_F\cap\bar\Omega$ with compact support which belong to the domain of the {\it actual} adjoint $\bar\partial^*$ of $\bar\partial$.  We say that the {\it order of subellipticity} for the scalar multiplier $F$ is $\geq\varepsilon$.  We are only interested in an effective lower bound for the order of subellipticity $\varepsilon_F$ and will not study the supremum of all possible such $\varepsilon_F$.

\medbreak The collection of smooth function germs at $P$ which are scalar multipliers at $P$ forms an ideal.  This ideal is called
the {\it multiplier ideal} at $P$ and is denoted by $I_P$.  The precise definition of the multiplier ideal makes precise the intuitive motivation that the direction and the order of vanishing
of the multiplier ideal measures the location and extent of the failure of the {\it a priori} subelliptic estimate.

\medbreak The test function $\varphi$ which is multiplied by the scalar multiplier $F$ to yield the modified subelliptic estimate
(1.4.1) is not a scalar and is a $(0,1)$-form with $n$ components (of which the normal component is $0$).  It is possible to get more information by using {\it vector multipliers}
instead of just scalar multipliers.  A smooth germ of $(1,0)$-form $\theta$ at a point of $\partial\Omega$ (defined on some open neighborhood $U_\theta$ of $P$ in ${\mathbb C}^n$)
is called a {\it vector multiplier} if
for some positive number $\varepsilon_\theta>0$ and some positive constant $C_{{}_\theta}$ the subelliptic estimate
$$
\||\theta\cdot\varphi|\|_{\varepsilon_{{}_\theta}}^2\leq C_{{}_\theta}\left(\|\bar\partial
\varphi\|^2+\|\bar\partial^*\varphi\|^2+\|\varphi\|^2\right)
$$
of order $\varepsilon_\theta$, {\it modified} by the dot product with $\theta$, holds for
for all smooth $(0,1)$-forms $\varphi$ on $U_\theta\cap\bar\Omega$ with compact support which belong to
the domain of the {\it actual} adjoint $\bar\partial^*$ of $\bar\partial$.  Here the dot product $\theta\cdot\varphi$ means the
pointwise inner product $\left<g,\,\bar\theta\right>$ of the two $(0,1)$-forms $\bar\theta$ and $g$ with respect to the usual Euclidean Hermitian inner product of ${\mathbb C}^n$.  In other words,
if $\theta=\sum_{j=1}^n\theta_j dz_j$ and $\varphi=\sum_{j=1}^n\varphi_{\bar j}d\bar z_j$ (where $z_1,\cdots,z_n$ are the global coordinates of ${\mathbb C}^n$), then $\theta\cdot \varphi=\sum_{j=1}^n\theta_j \varphi_{\bar j}$.  The convention to introduce a vector multiplier $\theta$ as a $(1,0)$-form and the dot product $\theta\cdot\varphi$, instead of introducing a $(0,1)$-form $\psi$ and the pointwise inner product $\left<\varphi,\,\psi\right>$, is chosen so that, in the case of a special domain described in (2.8.1) below, one needs only consider scalar multipliers and vector multipliers which are holomorphic (see (2.9.3) below).

\medbreak We say that the {\it order of subellipticity} for the vector multiplier $\theta$ is $\geq\varepsilon_{{}_\theta}$.  Again we are only interested in an effective lower bound for the order of subellipticity $\varepsilon_{{}_\theta}$ and will not study the supremum of all possible such $\varepsilon_{{}_\theta}$.

\medbreak The collection of smooth germs of $(1,0)$-forms at $P$ which are vector multipliers at $P$ forms a module over the algebra of all smooth function germs at $P$.  This module is called
the {\it module of vector multipliers} at $P$ and is denoted by $A_P$.

\bigbreak\noindent{\bf\S 2.} {\sc Generation of Multipliers in Kohn's Algorithm}

\bigbreak It is easy to define multipliers to describe the location and the extent of failure of subelliptic estimates, but it is difficulty to use multipliers to find easily verifiable conditions on the weakly pseudoconvex domain $\Omega$ to obtain subelliptic estimates and thereby the smoothness of the Kohn solution $u$ on $\bar\Omega$ from the smoothness of the right-hand side $f$ on $\bar\Omega$.  The most important part of the theory of multipliers for the complex Neumann problem is the generation of scalar and vector multipliers by Kohn's algorithm which makes it possible to study subelliptic estimates from the geometric condition (known as {\it finite type}) of the finiteness of the maximum normalized order of contact between the boundary $\partial\Omega$ of $\Omega$ and a local holomorphic curve $f:\Delta\to{\mathbb C}^n$ (where $\Delta$ is the open unit disk in ${\mathbb C}$).  Here the normalized order of contact means that the vanishing order of the pullback by $f$ of the defining function $r$ of $\partial\Omega$ divided by the vanishing order of $f$.  The precise definition of finite type condition and its history will be given later in (2.3) below.

\bigbreak\noindent(2.1) {\it Kohn's Algorithm to Generate Scalar and Vector Multipliers.} The following three procedures constitute the Kohn algorithm of generating scalar and vector multipliers at a boundary point $P$ of $\partial\Omega$.

\medbreak\noindent(A) {\small{\sc Initial Multipliers}}

\begin{itemize}\item[(i)] The function $r$ belongs to the ideal of multipliers $I_P$ at $P$. Its order of subellipticity for the scalar multiplier $r$ is $\geq 1$.

\item[(ii)] For any germ of smooth vector field $\xi=\sum_{k=1}^n a_k\frac{\partial}{\partial z_k}$ at $P$ of type $(1,0)$ which is tangential to $\partial\Omega$, the $(1,0)$-form
    $$(\partial\bar\partial r)\intprod\bar\xi=\sum_{j,k=1}^n\overline{\xi_k}\frac{\partial^2 r}{\partial z_j\partial\bar z_k}\,dz_j$$
    (which is the interior product of the $(1,1)$-form $\partial\bar\partial r$ and the $(0,1)$-vector field $\bar\xi$) belongs to
the module $A_P$ of vector multipliers at
$P$.  The order of subellipticity for the vector multiplier $\partial\bar\partial_j r$ is $\geq\frac{1}{2}$ for $1\leq j\leq n-1$ at points $P$ where $\partial r$ is normalized to
be $dz_n$.
(See [Kohn1979, p.97, (4.29)]).
\end{itemize}

\medbreak\noindent(B) {\small{\sc Generation of New Multipliers}}

\begin{itemize}\item[(i)]  If $f\in I_P$, then $\partial f\in A_P$.  The order $\varepsilon_{{}_{\partial f}}$ of subellipticity for the vector multiplier $\partial f$ is $\geq\frac{\varepsilon_f}{2}$ if the order of subellipticity of $f$ is $\geq\varepsilon_f$.  (See [Kohn1979, p.99, (4.42)].)

\item[(ii)] If $\theta_1,\cdots,\theta_{n-1}\in
A_P$ and $$\theta_1\wedge\cdots\wedge\theta_{n-1}\wedge\partial r= f_{\theta_1,\cdots,\theta_{n-1}}\,dz_1\wedge\cdots\wedge dz_n,$$ then $f_{\theta_1,\cdots,\theta_{n-1}}\in I_P$. In other words, for vector multipliers $\theta_1,\cdots,\theta_{n-1}$ the coefficient of $dz_1\wedge\cdots\wedge dz_n$ in $\theta_1\wedge\cdots\wedge\theta_{n-1}\wedge\partial r$ is a scalar multiplier.  Moreover, if the order of subellipticity of each $\theta_j$ (for $1\leq j\leq n-1$) is $\geq\varepsilon$, then the order of subellipticity of $f_{\theta_1,\cdots,\theta_{n-1}}$ is $\geq\varepsilon$.  We denote the scalar multiplier $f_{\theta_1,\cdots,\theta_{n-1}}$ by $\det_{n-1}(\theta_1,\cdots,\theta_{n-1})$.
\end{itemize}

\medbreak\noindent(C) {\small{\sc Real Radical Property}}

\begin{itemize}\item[] If $g\in I_P$ and
$\left|f\right|^m\leq\left|g\right|$ for some $m\in{\mathbb N}$, then $f\in I_P$.  The order of subellipticity for the scalar multiplier $f$ is $\geq\frac{\varepsilon}{m}$ if the order of subellipticity of $g$ is $\geq\varepsilon$.  (See [Kohn1979, p.98, Lemma 4.34].)
\end{itemize}

\bigbreak\noindent For the purpose of discussing the effectiveness of Kohn's algorithm later, we introduce now two terms concerning the radicals of ideals.  For an ideal $J$ of smooth function germs at $P$, we call the ideal of all smooth functions $f$ such that $|f|^m\leq|g|$ for some $m\in{\mathbb N}$ and for some $g\in J$ the {\it full real radical} of $J$.  If $q\in{\mathbb N}$ is given, we call the ideal of all smooth functions $f$ such that $|f|^q\leq|g|$ for some $g\in J$ the {\it real radical of root order $q$}.

\vskip.2in
\bigbreak\noindent(2.2) {\it Key Features of Kohn's Algorithm.}  The subelliptic estimate holds at a point $P$ of $\partial\Omega$ when there is
a scalar multiplier which is nonzero at $P$.  Kohn's algorithm seeks to reduce the vanishing order of scalar multipliers by differentiation and by
root-taking.  The procedure of differentiation described in (B) in (2.1) above would only allow certain differential operators to
lower the vanishing order of multipliers, namely only $(1,0)$-differentiation is allowed and only the determinants of coefficients of $(1,0)$-differentials
of scalar multipliers (from Cramer's rule) can be used to produce new scalar multipliers.  The procedure of root-taking described in (C) in (2.1) above identifies a smooth function germ as a scalar multiplier when a positive integral power of its absolute-value is dominated by the absolute-value of some known scalar multiplier.

\bigbreak\noindent(2.3) {\it Condition of Finite Type.}
The {\it type} $m$ at a point $P$ of the
boundary of weakly pseudoconvex $\Omega$ is the supremum of the
normalized touching order
$$\frac{{\rm ord}_0\left(r\circ f\right)}{{\rm ord}_0f},$$
to $\partial\Omega$, of all local holomorphic curves
$f:\Delta\to{\mathbb C}^n$ with $\varphi(0)=P$,
where $\Delta$
is the open unit $1$-disk and ${\rm ord}_0$ is the vanishing order
at the origin $0$.  The notion of finite type was first introduced by Kohn in 1972 [Kohn1972, p.525, Def.2.3] for the case of $n=2$ where the formulation is in terms of
the nonvanishing of $\partial r$ on the iterated Lie brackets of tangential vector fields to $\partial\Omega$ of types $(1,0)$ and $(0,1)$.  It was extended to
the case of a general $n$ by D'Angelo in [D'Angelo1979, p.59] in terms of finite algebraic obstructions to the existence of a nontrivial holomorphic complex curve germ in ${\mathbb C}^n$ tangential to $\partial\Omega$ and then in the formulation in terms of normalized touching order in [D'Angelo1982 , p.625, Definition 2.16 and Definition 2.18].

\bigbreak\noindent(2.4) {\it Kohn's Conjecture.}  The goal of the theory of multipliers for the complex Neumann problem is to prove by using the procedures of generating new multipliers to prove that the function identically $1$ is generated as a multiplier so that for a bounded weakly pseudoconvex domain $\Omega$ with finite type $m$ the subelliptic estimate for some positive order $\varepsilon>0$ holds and as a consequence the Kohn solution $u$ of the $\bar\partial$-equation $\bar\partial u=f$ on $\Omega$ with the right-hand side $f$ smooth on $\bar\Omega$ is also smooth on $\bar\Omega$.  Moreover, there is effectiveness in the use of the procedures to generate new multipliers so that $\varepsilon$ is some explicit function of the type $m$ of the domain $\Omega$ and its complex dimension $n$.

\bigbreak\noindent(2.5) {\it Full-Real-Radical Kohn Algorithm.}  If the question of effectiveness is to be set aside, the following algorithm can be used for the generation of new multipliers in $I_P$.

\medbreak\noindent(i) We start with the initial member $r$ of $I_P$ and denote by $I_P^{(0)}$ the ideal generated by this initial scalar multiplier $r$.  Likewise, we start out with the initial members $(\partial\bar\partial r)\intprod\bar\xi$ of $A_P$ for all choices of smooth $(1,0)$-vector fields $\xi$ tangential to $\partial\Omega$ and denote by $A_P^{(0)}$ the module generated by these initial vector multipliers.

\medbreak\noindent(ii) We use induction on the nonnegative integer $\nu$ to define $I_P^{(\nu)}$ and $A_P^{(\nu)}$ as follows.  We add to $I_P^{(\nu)}$ all the multipliers
    $\det_{n-1}(\theta_1,\cdots,\theta_{n-1})$ for $\theta_1,\cdots,\theta_{n-1}$ in $A_P^{(\nu)}$ to form an ideal which we denote by $\hat I_P^{(\nu)}$.  For the induction step from $\nu$ to $\nu+1$ we define $I_P^{(\nu+1)}$ as the {\it full real radical} of the ideal $\hat I_P^{(\nu)}$.  We add to $A_P^{(\nu)}$ all the vector multipliers $\partial F$ for $F\in I_P^{(\nu+1)}$ to form $A_P^{(\nu+1)}$.

\bigbreak We then hope that if the pseudoconvex domain $\Omega$ is of finite type $m$, the above construction of $I_P^{(\nu)}$ and $A_P^{(\nu)}$ by induction on $\nu$ will result in $I_P^{(\nu_m)}$ containing the function $1$ for some $\nu_m$ which effectively depends on $m$ and $n$.  In other words, the algorithm {\it terminates}.  In particular, the subelliptic estimate of some positive order $\varepsilon>0$ holds.  However, since there is no control on the order of the root-taking used in going from the ideal $\hat I_P^{(\nu)}$ to its {\it full real radical} $I_P^{(\nu+1)}$ there is no way for us to conclude that the order $\varepsilon$ of the subellipticity proved depends effectively on $m$ and $n$.  For the purpose of our discussion we call this algorithm of using induction on $\nu$ to construct $I_P^{(\nu)}$ with the goal of ending up with $1\in I_P^{(\nu_m)}$ the {\it full-real-radical Kohn algorithm}.

\bigbreak\noindent(2.6) {\it Effective Kohn Algorithm.}  In order to end up with an effective order $\varepsilon$ of subellipticity, it will turn out that one needs to modify the full-real-radical Kohn algorithm to use the following {\it effective Kohn algorithm}.

\medbreak\noindent(i) The starting point for the effective Kohn algorithm is the same as for the full-real-radical Kohn algorithm.  Again we start with the initial member $r$ of $I_P$ and denote by $I_P^{(0)}$ the ideal generated by this initial scalar multiplier $r$.  Likewise, we start out with the initial members $(\partial\bar\partial r)\intprod\bar\xi$ of $A_P$ for all choices of smooth $(1,0)$-vector fields $\xi$ tangential to $\partial\Omega$ and denote by $A_P^{(0)}$ the module generated by these initial vector multipliers.

\medbreak\noindent(ii) What is different with the effective Kohn algorithm is that we introduce a positive integer $q_\nu$ for every nonnegative integer index $\nu$ so that we take the real radical of root order $q_\nu$ instead of the full real radical in the step of going from $\nu-1$ to $\nu$.  First we set $q_0=1$ and set $I_P^{(0,1)}=I_P^{(0)}$.  This is motivated by the fact that the full real radical of $I_P^{(0)}$ is $I_P^{(0)}$ again because $dr$ is nowhere zero on $\partial\Omega$.  We use induction on the nonnegative integer $\nu$ to define a sequence of positive integers $q_\nu$ and then $I_P^{(\nu,q_\nu)}$ and $A_P^{(\nu,q_\nu)}$ as follows.  We add to $I_P^{(\nu,q_\nu)}$ all the multipliers
    $\det_{n-1}(\theta_1,\cdots,\theta_{n-1})$ for $\theta_1,\cdots,\theta_{n-1}$ in $A_P^{(\nu,q_\nu)}$ and then use the resulting collection of multipliers to form an ideal which we denote by $\hat I_P^{(\nu,q_\nu)}$.  For the induction step from $\nu$ to $\nu+1$ we define $I_P^{(\nu+1,q_{\nu+1})}$ as the {\it real radical} of root order $q_{\nu+1}$ of the ideal $\hat I_P^{(\nu,q_\nu)}$.  We add to $A_P^{(\nu,q_\nu)}$ all the vector multipliers $\partial F$ for $F\in I_P^{(\nu+1,q_{\nu+1})}$ and then use the resulting collection of vector multipliers to form a module which is $A_P^{(\nu+1,q_{\nu+1})}$ for the next step of going from $\nu$ to $\nu+1$ in the construction by induction on $\nu$.

\bigbreak We then hope that if the pseudoconvex domain $\Omega$ is of finite type $m$, the above construction of $I_P^{(\nu,q_\nu)}$ and $A_P^{(\nu,q_\nu)}$ by induction on $\nu$ will result in $I_P^{(\nu_m,q_{\nu_m})}$ containing the function $1$ for some $\nu_m$ which effectively depends on $m$ and $n$.  Moreover, we hope that each $q_\nu$ for $0\leq\nu\leq\nu_m$ depends also effectively on $m$ and $n$ so that the subelliptic estimate of some positive order $\varepsilon>0$ holds with $\varepsilon$ depending effectively on $m$ and $n$.

\medbreak Note that $I_P^{(\nu,\mu)}$ and $A_P^{(\nu,\mu)}$ are defined only when $\mu=q_\nu$.  Instead of using the notations $I_P^{(\nu,q_\nu)}$ and $A_P^{(\nu,q_\nu)}$ we could have used the notations which depend only on $\nu$, for example, $\tilde I_P^{(\nu)}$ and $\tilde A_P^{(\nu)}$.  We prefer the clumsier notations $I_P^{(\nu,q_\nu)}$ and $A_P^{(\nu,q_\nu)}$ to highlight the effective choice of $q_\nu$.

\medbreak In the effective Kohn algorithm described above, no procedure is given to determine the sequence $q_\nu$.  The sequence $q_\nu$ is obtained by a rather complicated algebraic geometric argument.  Since the purpose is to achieve $1\in I_P^{(\nu_m,q_{\nu_m})}$, it suffices to describe explicitly how to use algebraic geometric techniques to determine the procedures of differentiation and root-taking of order $\leq q_\nu$ to construct scalar multipliers and vector multipliers from some initial scalar multipliers and vector multipliers to achieve $1\in I_P^{(\nu_m,q_{\nu_m})}$.

\medbreak Kohn's papers [Kohn1977, Kohn1979] discussed relations between subelliptic estimates, finite type property, and the termination of the full-real-radical Kohn algorithm for smooth weakly pseudoconvex domains.  The relations can be summarized in the following Kohn conjecture, formulated separately for the full-real-radical Kohn algorithm and the effective Kohn algorithm.

\bigbreak\noindent(2.7) {\it Conjecture on Full-Real-Radical Kohn Algorithm.} For a bounded weakly pseudoconvex domain $\Omega$ in ${\mathbb C}^n$ with smooth boundary and of finite type $m$ and for a point $P$ of the boundary of $\Omega$, in the ascending chain of multiplier ideals $I_P^{(\nu)}$ for $\nu\in{\mathbb N}\cup\{0\}$ there exists some $\nu^*\in{\mathbb N}$ such that $I_P^{(\nu^*)}$ contains the constant function $1$.  In other words, the chain of multiplier ideals $I_P^{(\nu)}$ for $\nu\in{\mathbb N}\cup\{0\}$ terminates at $\nu=\nu^*$.

\medbreak If the positive integer $\nu^*$ depends effectively on $m$ and $n$, then we say that the full-real-radical Kohn algorithm terminates effectively.

\bigbreak\noindent(2.8) {\it Conjecture on Effective Kohn Algorithm.} For a bounded weakly pseudoconvex domain $\Omega$ in ${\mathbb C}^n$ with smooth boundary and of finite type $m$ and for a point $P$ of the boundary of $\Omega$, there exist $\tilde\nu\in{\mathbb N}$ and a sequence of positive numbers $q_1,\cdots,q_{\tilde\nu}$ such that in the ascending chain of multiplier ideals $I_P^{(\nu,q_\nu)}$ for $0\leq\nu\leq\tilde\nu$, the multiplier ideal $I_P^{(\tilde\nu,q_{\tilde\nu})}$ contains the constant function $1$.  Moreover, the positive integer $\tilde \nu$ and the sequence of positive integers $q_1,\cdots,q_{\tilde\nu}$ depend effectively on $m$ and $n$.

\bigbreak For the conjecture on the full-real-radical Kohn algorithm, when the boundary $\partial\Omega$ of the bounded weakly pseudoconvex domain $\Omega$ is assumed real-analytic, the finite type condition becomes a conclusion instead of an assumption.  Kohn first showed that if the the chain of multiplier ideals $I_P^{(\nu)}$ for $\nu\in{\mathbb N}\cup\{0\}$ does not terminate, then the boundary $\partial\Omega$ contains a local real-analytic subvariety of holomorphic dimension $\geq 1$ ([Kohn1977, p.2215, Lemma 20] and [Kohn1979, p.113, Proposition 6.20]).  Then Diederich-Fornaess proved that the real-analytic boundary of a bounded weakly pseudoconvex domain cannot contain a local real-analytic subvariety of holomorphic dimension $\geq 1$ ([Diederich-Fornaess1978, p.373, Lemma 2] and [Diederich-Fornaess1978, p.374, Theorem 3]). This result of Kohn and Diederich-Fornaess holds not only for the case of $(0,1)$-forms but for the general case of $(0,q)$-forms for $1\leq q\leq n$.  Since their method of proof is by contradiction, there is no effectiveness in the termination of the chain of multiplier ideals $I_P^{(\nu)}$.

\bigbreak\noindent(2.8.1) For the conjecture on Kohn's effective algorithm, the paper [Siu2010] introduces algebraic geometric techniques to study the problem by looking first at the special case of $\Omega$ being a special domain.  For notational convenience we now consider a domain $\Omega$ in ${\mathbb C}^{n+1}$ instead of ${\mathbb C}^n$.  A {\it special domain} $\Omega$ in
${\mathbb C}^{n+1}$ (with coordinates $w,z_1,\cdots,z_n$) is a
bounded domain given by
$$
{\rm
Re\,}w+\sum_{j=1}^N\left|F_j\left(z_1,\cdots,z_n\right)\right|^2<0,\leqno{(2.8.1.1)}
$$
where $F_j\left(z_1,\cdots,z_n\right)$ which is defined on some open
neighborhood of $\bar\Omega$ in ${\mathbb C}^{n+1}$ depends only on
the variables $z_1,\cdots,z_n$ and is holomorphic in
$z_1,\cdots,z_n$ for each $1\leq j\leq N$.   In [Siu2010] the verification of the conjecture on Kohn's effective algorithm for the case of $n=2$ (which means for special domains of complex dimension $3$) was given in detail, with only indications for the case of special domains of general complex dimension.  A rough outline was given there for the extension of the method of algebraic geometric techniques first to the general real-analytic case and then the general smooth case.  We carry out here the effective Kohn algorithm, which was introduced in [Siu2010], in a way which keeps track of the order of subellipticity in each step, in the context of comparing the full-real-radical Kohn algorithm and the effective Kohn algorithm.

\bigbreak In the chain of multiplier ideals $I_P^{(\nu)}$ in the conjecture for the full-real-radical Kohn algorithm, if for every $\nu$ there is an effective positive integer $p_\nu$ ({\it i.e.,} dependent only on $m$ and $n$) such that $\left(I_P^{(\nu+1)}\right)^{p_\nu}\subset I_P^{(\nu)}$, then we can use $p_\nu=q_\nu$ and the proof of the conjecture on Kohn's effective algorithm is simply reduced to the conjecture on the full-real-radical Kohn algorithm with effective termination.  Unfortunately, there are simple examples, even for special domains of complex dimension $3$, with no effective $p_\nu$. This means that the conjecture for effective Kohn algorithm is different from the conjecture for the full-real-radical Kohn algorithm with effective termination.  A simple example of this kind was given by Catlin and D'Angelo in [Catlin-D'Angelo2010], which we will discuss in \S4 below.

\bigbreak\noindent(2.9) {\it Algebraic Geometric Techniques in Effective Kohn Algorithm.}  We now explain how the algebraic geometric techniques precisely work to provide naturally the positive integer $q_\nu$ to make the process effective. The reason for considering special domains in ${\mathbb C}^{n+1}$ is that instead of ideals of smooth function germs on ${\mathbb C}^{n+1}$ we need only consider ideals of holomorphic function germs on ${\mathbb C}^n$.  That is the reason why for notational convenience we suddenly consider domains in ${\mathbb C}^{n+1}$ instead of domains in ${\mathbb C}^n$.  We focus on the case where each $F_j$ vanish at the origin of ${\mathbb C}^n$ (for $1\leq j\leq N$) so that the origin belongs to the boundary of $\Omega$.  We are concerned only with the problem of the subelliptic estimate at the origin.

\medbreak The new notion of pre-multipliers needs to be introduced.  A holomorphic function germ $f(z_1,\cdots,z_n)$ on ${\mathbb C}^n$ at the origin is a {\it pre-multiplier} if its differential $df$ is a vector multiplier at the origin.  We now verify that the holomorphic function germs $F_1,\cdots,F_N$ are pre-multipliers and the order of subellipticity of each $dF_j$ is $\geq\frac{1}{4}$.

\bigbreak\noindent(2.9.1) {\it Levi Form and Initial Vector Multiplier for Special Domain.}  First, for a special domain we write down explicit expressions for a tangent vector, its Levi form, and a smooth test $(1,0)$-form in the domain of the actual adjoint $\bar\partial^*$ of $\bar\partial$.  Since the defining function for the special domain $\Omega$ is
$$
r={\rm
Re\,}w+\sum_{j=1}^N\left|F_j\left(z_1,\cdots,z_n\right)\right|^2,
$$
it follows that
$$\partial r=\frac{1}{2}\,dw+\sum_{j=1}^N \overline{F_j}\,dF_j
$$
and all the $(1,0)$-vector tangential to $\partial\Omega$ at the origin are of the form
$$b\frac{\partial}{\partial w}+\sum_{k=1}^n a_j\frac{\partial}{\partial z_k}$$ for $b,a_1,\cdots,a_n\in{\mathbb C}$ with
$$
b=-2\sum_{1\leq j\leq N,\,1\leq k\leq n}a_k\,\overline{F_j}\,\frac{\partial F_j}{\partial z_k}.
$$
The value of the Levi-form $\partial\bar\partial r=\sum_{j=1}^N dF_j\wedge\overline{dF_j}$ at $\eta\wedge\bar\eta$ for the element
$$
\eta=-2\left(\sum_{1\leq j\leq N,\,1\leq k\leq n}a_k\,\overline{F_j}\,\frac{\partial F_j}{\partial z_k}\right)\frac{\partial}{\partial w}+\sum_{k=1}^n a_j\frac{\partial}{\partial z_k}
$$
of
$T^{(1,0)}_{\partial\Omega}$ is equal to
$$
\sum_{1\leq j\leq N,\,1\leq k\leq n}\left|a_k\,\frac{\partial F_j}{\partial z_k}\right|^2,
$$
which means that a special domain is always weakly pseudoconvex and is strict pseudoconvex at a boundary point if and only if $dF_{j_1},\cdots,dF_{j_n}$ are ${\mathbb C}$-linearly independent at that point for some $1\leq j_1<\cdots<j_n\leq N$.
For the element
$$\xi=-2\left(\sum_{j=1}^N\overline{F_j}\,\frac{\partial F_j}{\partial z_\mu}\right)\frac{\partial}{\partial w}+\frac{\partial}{\partial z_\mu}
$$
of $T^{(1,0)}_{\partial\Omega}$ (for some fixed $1\leq\mu\leq n$) the $(1,0)$-form
$$(\partial\bar\partial r)\intprod\bar\xi=\sum_{\nu=1}^n\left(\sum_{j=1}^N(\partial_{\nu}F_j)\overline{(\partial_\mu F_j)}\right)dz_\nu$$
is a vector multiplier.

\medbreak For an open neighborhood $U$ of a point of $\partial\Omega$ in ${\mathbb C}^n$, if $\varphi=\sum_{\nu=1}^n\varphi_{\bar\nu}d\bar z_\nu+\hat\varphi d\bar w$ is a smooth test $(0,1)$-form
on $U\cap\bar\Omega$ with compact support which belongs to the
domain of the {\it actual} adjoint $\bar\partial^*$ of $\bar\partial$, then the normal component of $\varphi$ vanishes on $U\cap\partial\Omega$, which means that
$$
\begin{aligned}
(\partial r)\cdot\varphi&=\left(\frac{1}{2}\,dw+\sum_{j=1}^N \overline{F_j}\,dF_j\right)\cdot\left(\sum_{\nu=1}^n\varphi_{\bar\nu}d\bar z_\nu+\hat\varphi d\bar w\right)\cr
&=\frac{1}{2}\,\hat\varphi+\sum_{j=1}^N\sum_{\nu=1}^n\overline{F_j}\left(\partial_\nu F_j\right)\varphi_{\bar\nu}\cr
\end{aligned}
$$
vanishes on $U\cap\partial\Omega$, or
$$
\hat\varphi=-2\sum_{j=1}^N\sum_{\nu=1}^n\overline{F_j}\left(\partial_\nu F_j\right)\varphi_{\bar\nu}
$$
on $U\cap\partial\Omega$.

\medbreak It will turn out that in the case of a special domain the use of multipliers and vector multipliers can be limited to those which are holomorphic in $z_1,\cdots,z_n$ and are independent of $w$ and as a consequence in the study of the subelliptic estimate for a special domain the component $\hat\varphi$ of a test $(0,1)$-form $\varphi$ actually plays no role.  (See (2.9.3) below.)

\bigbreak\noindent(2.9.2) {\it Initial Pre-Multiplier for Special Domain.}
To verify that $dF_j$ is a vector multiplier with order of subellipticity $\geq\frac{1}{4}$, we take an open neighborhood $U$ of $0$ on ${\mathbb C}^{n+1}$ on which the holomorphic functions $F_1,\cdots,F_N$ of $(z_1,\cdots,z_n)\in{\mathbb C}^n$ are defined. For $0<\varepsilon\leq\frac{1}{2}$ and any smooth $(0,1)$-form $\varphi=\sum_{\nu=1}^n\varphi_{\bar\nu}d\bar z_\nu+\hat\varphi d\bar w$ on $U\cap\bar\Omega$ with compact support which belongs to the
domain of the {\it actual} adjoint $\bar\partial^*$ of $\bar\partial$,
we have
$$
\begin{aligned}
&\sum_{j=1}^N\left\|\left|\sum_{\nu=1}^n(\partial_\nu F_j)\varphi_{\bar\nu}\right|\right\|_\varepsilon^2
\cr
&=\sum_{j=1}^N\sum_{\nu,\mu=1}^n\left(\Lambda^\varepsilon\left((\partial_\nu F_j)\varphi_{\bar\nu}\right),\,\Lambda^\varepsilon\left((\partial_\mu F_j)\varphi_{\bar\mu}\right)\right)
\cr
&=\sum_{j=1}^N\sum_{\nu,\mu=1}^n\left(\overline{(\partial_\mu F_j)}
\Lambda^{2\varepsilon}\left((\partial_\nu F_j)\varphi_{\bar\nu}\right),\,\varphi_{\bar\mu}\right)
\cr
&=\sum_{j=1}^N\sum_{\nu,\mu=1}^n\left(\Lambda^{2\varepsilon}\left(\overline{(\partial_\mu F_j)}
(\partial_\nu F_j)\varphi_{\bar\nu}\right),\,\varphi_{\bar\mu}\right)+\sum_{j=1}^N\sum_{\nu,\mu=1}^n\left(\left[\overline{(\partial_\mu F_j)},\Lambda^{2\varepsilon}\right]\left(
(\partial_\nu F_j)\varphi_{\bar\nu}\right),\,\varphi_{\bar\mu}\right)
\cr
&=\sum_{\mu=1}^n\left(\Lambda^{2\varepsilon}\left(\sum_{\mu=1}^n\left(\left(\sum_{j=1}^N\overline{(\partial_\mu F_j)}
(\partial_\mu F_j)\right)\varphi_{\bar\nu}\right)\right),\,\varphi_{\bar\mu}\right)\cr
&\quad+\sum_{j=1}^N\sum_{\nu,\mu=1}^n\left(\left[\overline{(\partial_\mu F_j)},\Lambda^{2\varepsilon}\right]\left(
(\partial_\nu F_j)\varphi_{\bar\nu}\right),\,\varphi_{\bar\mu}\right)
\cr
&\leq\left(\sum_{\mu=1}^n\left\|\left|\sum_{\nu=1}^n\left(\sum_{j=1}^N(\partial_\nu F_j)\overline{(\partial_\mu F_j)}\right)\varphi_{\bar\nu}\right|\right\|_{2\varepsilon}^2\right)^{\frac{1}{2}}\left\|\varphi\right\|+C_1\left\|\varphi\right\|^2\cr
&\leq C_2\left(\|\bar\partial
\varphi\|^2+\|\bar\partial^*\varphi\|^2+\|\varphi\|^2\right)\cr
\end{aligned}
$$
where $C_1$ and $C_2$ are positive constants independent of $\varphi$ (but depend on $U$ and $\varepsilon$), because
$$
\sum_{\nu=1}^n\left(\sum_{j=1}^N(\partial_{\nu}F_j)\overline{(\partial_\mu F_j)}\right)dz_\nu
$$
is a vector multiplier whose order of subellipticity is $\geq\frac{1}{2}$ for $1\leq\nu\leq n$.  This finishes the verification that each $dF_j=\sum_{\nu=1}^n(\partial_\nu F_j)dz_j$ is a vector multiplier whose order of subellitpicity is $\geq\frac{1}{4}$.

\bigbreak\noindent(2.9.3) {\it Holomorphic Multipliers for Special Domain.} We start out with pre-multipliers $F_1,\cdots,F_N$ and the vector multipliers $dF_1,\cdots,dF_N$ which are holomorphic $1$-forms in the variables $z_1,\cdots,z_n$ obtained from them.  For the special domain $\Omega$ we will only work with vector multipliers which are holomorphic $1$-forms in the variables $z_1,\cdots,z_n$.  When we have $n$ such vector multipliers
$G^{(j)}=\sum_{\nu=1}^n G^{(j)}(z_1,\cdots,z_n)dz_j$ of order of subellipticity $\geq\varepsilon_\nu$ for $1\leq\nu\leq n$, when we apply the procedure (B)(ii) in (2.1) above to generate new multipliers, from
$$\partial r=\frac{1}{2}\,dw+\sum_{j=1}^N \overline{F_j}\,dF_j$$
we get
$$
G^{(1)}\wedge\cdots\wedge G^{(n)}\wedge\partial r=\frac{1}{2}\det\left(G^{(\nu)}_j\right)_{1\leq\nu,j\leq n}dz_1\wedge\cdots\wedge dz_n\wedge dw
$$
to conclude that the holomorphic function
$$
\det\left(G^{(\nu)}_j\right)_{1\leq\nu,j\leq n}
$$
of the variables $z_1,\cdots,z_n$ is a scalar multiplier whose order of subellipticity is $\geq\min(\varepsilon_1,\cdots,\varepsilon_n)$.  So for the special domain $\Omega$, when we start out only with vector multipliers $dF_1,\cdots,dF_N$ (whose orders of subellipticity are all $\geq\frac{1}{4}$) and use only the procedures in B(i)(ii) and C described in (2.1) above to generate new scalar and vector multipliers, we need only work with scalar and vector multipliers which are holomorphic functions or holomorphic $1$-forms in the variables $z_1,\cdots,z_n$.

\medbreak We now translate the algorithm from the language of analysis to the language of algebraic geometry.  The procedures in the algorithm now read as follows.

\bigbreak\noindent(2.9.4) {\it Algebraic Geometric Formulation of Kohn Algorithm for Special Domain.}
The multiplicity $q$ of the ideal generated by $F_1,\cdots,F_N$ given by
$$
\dim_{\mathbb C}\left({\mathcal O}_{{\mathbb C}^n,0}\left/\sum_{j=1}^N
{\mathcal O}_{{\mathbb C}^n,0} F_j\right.\right)
$$
is related to the type $m$ of the special domain $\Omega$ at the origin by
$2m\leq q\leq(n+2)2m$.  See [Siu2010, Lemma(I.3) and Lemma(I.4)].  In particular, the special domain is of finite type at $0$ if and only if $0$ is an isolated point of
the common zero-set of $F_1,\cdots,F_N$. Assume that the multiplicity $q$ of the ideal generated by $F_1,\cdots,F_N$ is finite.  The $N$ pre-multipliers $F_1,\cdots,F_N$ (with order of subellipticity of each $dF_j$ at least $\frac{1}{4}$) are all that we start out with.  For a special domain of complex dimension $n+1$ there are only the following two procedures from Kohn's algorithm.

\medbreak\noindent(i) If holomorphic function germs $g_1,\cdots,g_n$ at $0$ on ${\mathbb C}^n$ are pre-multipliers (which automatically include all multipliers), then the coefficient of $dz_1\wedge\cdots\wedge dz_n$ in $dg_1\wedge\cdots\wedge dg_n$ is a multiplier.  In other words, the Jacobian determinant
$$
\frac{\partial(g_1,\cdots,g_n)}{\partial(z_1,\cdots,z_n)}
$$
of the holomorphic functions $g_1,\cdots,g_n$ with respect to the variables $z_1,\cdots,z_n$ is a multiplier.
If the order of subellipticity of each $g_1,\cdots,g_n$ is $\geq\eta$, then the order of subellipticity of their Jacobian determinant is $\geq\frac{\eta}{2}$.

\medbreak\noindent(ii) If $g$ and $f$ are holomorphic function germs at $0$ on ${\mathbb C}^n$ with $f^m=g$ for some positive integer $m$ and if $g$ is a multiplier whose order of subellipticity is $\geq\eta$, then $f$ is also a multiplier whose order of subellipticity is $\geq\frac{\eta}{m}$.

\medbreak Note that the set of all multipliers forms an ideal in the ring of all holomorphic function germs and the set of all vector-multipliers form a module over the ring of all holomorphic function germs, but in general the set of all pre-multipliers does {\it not} form a module over the ring of all holomorphic functions, because though the differential $dF$ of a pre-multiplier $F$ is a multiplier, yet for any holomorphic function germ $g$ the differential $d(gF)$ is equal to $gdF+Fdg$ and the term $Fdg$, unlike the term $gdF$, is in general {\it not} a vector-multiplier.

\medbreak We now discuss the algebraic geometric formulations of the steps in both the full-real-radical Kohn algorithm and the effective Kohn algorithm.  We will describe first the steps in the full-real-radical Kohn algorithm.  Then we will describe the effective Kohn algorithm but only in the case when the special domain is of complex dimension $3$.

\bigbreak\noindent(2.10) {\it Steps in Full-Real-Radical Kohn Algorithm.}   We start out with the ${\mathbb C}$-vector space $V_0$ of initial pre-multipliers generated by $F_1,\cdots,F_N$.  Let $J_0$ be the ideal generated by all the Jacobian determinants of any $n$ elements $g_1,\cdots,g_n$ of $V_0$.  Let $I_0$ be the radical of $J_0$.

\medbreak Inductively we construct the ${\mathbb C}$-vector space $V_\nu$, the ideal $J_\nu$, and the ideal $I_\nu$ for any nonnegative integer $\nu$ as follows.  For the step of going from $\nu$ to $\nu+1$, we let $V_{\nu+1}$ be the ${\mathbb C}$-vector space generated by elements of $V_\nu$ and all elements of the ideal $I_\nu$.  Let $J_{\nu+1}$ be the ideal generated by all the Jacobian determinants of any $n$ elements of $V_{\nu+1}$.  Let $I_{\nu+1}$ be the radical of the ideal of $J_{\nu+1}$.  This finishes the construction by induction.

\medbreak Let $p_\nu$ be the smallest positive integer such that $I_\nu^{p_\nu}$ is contained in $J_\nu$.
Note that each element of $I_\nu$ is a multiplier, but each element of $V_\nu$ is only a pre-multiplier.

\medbreak The full-real-radical Kohn algorithm for $F_1,\cdots,F_N$ {\it terminates} if there exists some nonnegative $\tilde\nu$ such that $I_{\tilde\nu}$ is the unit ideal, which means the entire ring of all holomorphic function germs.  In that case we choose $\tilde\nu$ to be the smallest such nonnegative integer.

\medbreak We say that the full-real-radical Kohn algorithm for $F_1,\cdots,F_N$ terminates {\it effectively} if $\tilde\nu$ and each $p_\nu$ for $0\leq\nu\leq\tilde\nu$ are bounded by explicit functions of $n$ and $q$.

\medbreak The order of subellipticity for the multiplier $1$ is at least
$$
\frac{1}{2^{\tilde\nu+2}\prod_{\nu=0}^{\tilde\nu}p_\nu}.
$$
This order of subellipticity from the termination of the original Kohn algorithm is effective only when $\tilde\nu$ is effective and each $p_\nu$ is effective for each $0\leq\nu\leq\tilde\nu$.

\bigbreak\noindent(2.11) {\it Ideal Containing an Effective Power of its Radical in Effective Kohn Algorithm.}
The key difference between the full-real-radical Kohn algorithm (which in general is not effective) by the effective Kohn algorithm is to replace the taking of the radical $I_\nu$ of the ideal $J_\nu$ by an appropriately chosen sub-ideal $\tilde I_\nu$ of $I_\nu$ with the property that an effective power $(\tilde I_\nu)^{s_\nu}$ of $\tilde I_\nu$ is contained in $J_\nu$.   The choice of the sub-ideal $\tilde I_\nu$ of $I_\nu$ and the positive integer $s_\nu$ involves rather complicated algebraic geometric techniques.
In order to facilitate the explicit comparison of the full-real-radical Kohn algorithm with the effective Kohn algorithm in a concrete example (such as the example of Catlin-D'Angelo [Catlin-D'Angelo2010]), in the description of the steps in the effective Kohn algorithm for special domains we will confine ourselves to special domains of complex dimension $3$, which means the case of $n=2$.

\bigbreak\noindent{\bf\S 3.} {\sc Orders of Subellpiticity in Algebraic Geometric Techniques for 3-Dimensional Special Domain}

\bigbreak We now describe the steps in effective Kohn algorithm for special domains of complex dimension 3.  Again we start with holomorphic function germs $F_1,\cdots,F_N$ at $0$ on ${\mathbb C}^2$ (which define the special domain in ${\mathbb C}^3$ and which generates an ideal in ${\mathcal O}_{{\mathbb C}^2,0}$ of multiplicity $\leq q$).

\bigbreak\noindent(3.1) {\it Step One.}
We take two generic ${\mathbb C}$-linear combinations $\hat F_1, \hat F_2$ of $F_1,\cdots,F_N$ of vanishing order $\leq q$ at $0$ such that the Jacobian determinant
$$
h_2^*:=\frac{\partial(\hat F_1,\hat F_2)}{\partial(z_1,z_2)}
$$
has vanishing order $\leq q$ at $0$ as a holomorphic function germ at $0$ on ${\mathbb C}^2$.  An order of subellipticity of $h_2^*$ as a multiplier at the origin is at least $\frac{1}{4}$. Let
$$
h_2^*=\left(h_{2,1}^*\right)^{k_1}\cdots\left(h_{2,\ell_2}^*\right)^{k_{\ell_2}}
$$
be the factorization into irreducible holomorphic function germs $h_{2,1}^*,\cdots,h_{2,\ell_2}^*$ with $k_1\geq k_2\geq\cdots\geq k_{\ell_2}\geq 1$.  Since the vanishing order of $h_2^*$ is $\leq q$, we have $k_1+\cdots+k_{\ell_2}\leq q$ and, in particular, $k_1\leq q$. The holomorphic function germ
$$
\left(h_{2,1}^*\cdots h_{2,\ell_2}^*\right)^{k_1}
$$
contains $h_2^*$ as a factor and is therefore a multiplier whose order of subellipticity is $\geq\frac{1}{4}$.  Let
$$
\hat h_2=h_{2,1}^*\cdots h_{2,\ell_2}^*.
$$
Since the $k_1$-th power of $\hat h_2$ is a multiplier, it follows from $k_1\leq q$ that $\hat h_2$ is a multiplier whose order of subellipticity is $\geq\frac{1}{4q}$.  The construction of $\hat h_2$ from $h_2^*$ is to make sure that the divisor of $\hat h_2$ is reduced (which means that its multiplicity at any of its regular points is $1$, though it is possibly reducible with many branches).  Up to this point the only goal accomplished is to produce a {\it reduced} holomorphic function germ $\hat h_2$ at $0$, with vanishing order $\leq q$ at $0$, which is a multiplier whose order of subellipticity is $\geq\frac{1}{4q}$.  This step of construction of $\hat h_2$ from $h_2^*$ is included only for the sake of convenience and is not actually absolutely necessary.

\bigbreak\noindent(3.2) {\it Step Two.}  Take a generic ${\mathbb C}$-linear combination $h_1$ of $F_1,\cdots,F_N$ and a generic ${\mathbb C}$-linear coordinate system $w_1,w_2$ of ${\mathbb C}^2$ such that the following three conditions (i), (ii) and (iii) are satisfied.  Here generic means that $h_1=\sum_{j=1}^N a_jF_j$ and $w_k=\sum_{\ell=1}^2 g_{k\ell}z_\ell$ with the element
$(a_j, g_{k\ell})_{1\leq j\leq N,\,1\leq k,\ell\leq 2}$ of ${\mathbb C}^{N+1}$
chosen outside some proper subvariety ${\mathcal Z}$ of ${\mathbb C}^{N+4}$.  The proper subvariety ${\mathcal Z}$ of ${\mathbb C}^{N+4}$ can be obtained as the union of three proper subvarieties of ${\mathbb C}^{N+4}$, one for each of the three conditions (i), (ii) and (iii).

\medbreak\noindent(i) The origin $0$ of ${\mathbb C}^2$ is an isolated zero of $w_1$ and the holomorphic function germs $h_1$.

\medbreak\noindent(ii) The multiplicity at $0$ of the ideal generated by $h_1$ and $\hat h_2$ is $\leq q^2$.

\medbreak\noindent(iii) The multiplicity at $0$ of the ideal generated by $\hat h_2$ and $\left(\frac{\partial h_1}{\partial w_1}\right)_{w_2={\rm const}}$ is $\leq 3q^2$.

\bigbreak The reason why a generic ${\mathbb C}$-linear combination $h_1$ of $F_1,\cdots,F_N$ satisfies Condition (ii) is that the multiplicity of the ideal generated by $F_1,\cdots,F_N$ at the origin is $\leq q$.   The number $q^2$ in Condition (iii) is the product of the multiplicity $q$ of $\hat h_1$ at $0$ and the multiplicity $q$ of the ideal generated by $F_1,\cdots,F_N$ at $0$.

\bigbreak A choice of a generic ${\mathbb C}$-linear combination $h_1$ of $F_1,\cdots,F_N$ and a generic ${\mathbb C}$-linear coordinate system $w_1,w_2$ of ${\mathbb C}^2$ satisfying Condition (iii) is obtained from the following statement.

\medbreak\noindent(3.2.1) For each holomorphic function germ $f$ on ${\mathbb C}^n$ at $0$ which vanishes at $0$, the function germ $f^{n+1}$ on ${\mathbb C}^n$ at $0$ belongs to the ideal generated by $\frac{\partial f}{\partial z_j}$ for $1\leq j\leq n$ (where $z_1,\cdots,z_n$ are the coordinates of ${\mathbb C}^n$).

\medbreak\noindent The statement (3.2.1) is a consequence of Skoda's result on ideal generation [Skoda1972] (see [Siu2010, p.1232, Propositioni(A.2)]) and can be considered as a generalization, to general holomorphic function germs, of Euler's formula expressing a homogeneous polynomials in terms of its first-order partial derivatives.  We apply the statement to $f=h_1$ to conclude that $(h_1)^3$ belongs to the ideal generated by $\frac{\partial h_1}{\partial z_1},\,\frac{\partial h_1}{\partial z_2}$.  Condition (ii) implies that a generic ${\mathbb C}$-linear combination $H=\sum_{k=1}^2 c_k\,\frac{\partial h_1}{\partial z_k}$ satisfies the condition that the multiplicity at $0$ of the ideal generated by $H$ and $\hat h_2$ is $\leq 3q^2$.  We can choose a generic ${\mathbb C}$-linear coordinate system $w_1,w_2$ of ${\mathbb C}^2$ such that $H=\left(\frac{\partial h_1}{\partial w_1}\right)_{w_2={\rm const}}$.

\medbreak We are now ready to construct more multipliers from the choice of $h_1$ and $w_1,w_2$.  Since $h_1$ vanishes at $0$, it follows Condition (iii) that there exist holomorphic function germs $\alpha,\beta$ at $0$ on ${\mathbb C}^2$ such that
$$(h_1)^{3q^2}=\alpha\hat h_2+\beta\left(\frac{\partial h_1}{\partial w_1}\right)_{w_2={\rm const}}.
$$
Though $h_1$ which is a ${\mathbb C}$-linear combination $F_1,\cdots,F_N$ of pre-multipliers $F_1,\cdots,F_N$ may not be a multiplier, the long key argument given below is to show that actually $h_1$ is a multiplier.
We need the following statement concerning Weierstrass polynomials.

\bigbreak\noindent(3.2.2) {\it Weierstrass Polynomial for Image Curve Under Branched-Cover Map.}  Let $\zeta_1$ and $\zeta_2$  be holomorphic function germs at $0$ on ${\mathbb C}^2$ vanishing at $0$ such that the origin $0$ of ${\mathbb C}^2$ is an isolated point of the common zero-set of $\zeta_1$ and $\zeta_2$.  Let $H$ be a holomorphic function germ at $0$ on ${\mathbb C}^2$ vanishing at $0$ such that the origin $0$ of ${\mathbb C}^2$ is an isolated point of the common zero-set of $H$ and $\zeta_1$. Let $\ell$ be a positive integer.  If $\zeta_2^\ell$ belongs to the ideal generated by $H$ and $\zeta_1$, then some holomorphic function germ on ${\mathbb C}^2$ at $0$ of the form of a Weierstrass polynomial
$$
\zeta_2^\ell+\sum_{j=0}^{\ell-1}\theta_j(\zeta_1)\zeta_2^j
$$
(with $\zeta_1,\zeta_2$ as variables) contains $H$ as a factor, where $\theta_j$ is a holomorphic function germ on ${\mathbb C}$ at $0$ which vanishes at $0$ for $0\leq j\leq\ell-1$.

\medbreak For the proof of the statement (3.2.2) on Weierstrass polynomials, first we observe that for the special case where $\zeta_1,\zeta_2$ are the coordinate functions $z_1,z_2$ of ${\mathbb C}^2$ and the restriction of $H$ to $\{z_1=0\}$ is equal to $z_2^\ell$ as holomorphic function germ on ${\mathbb C}$ with $z_2$ as coordinate, the statement is simply the usual factorization of a holomorphic function germ $H$ as a product of a nowhere zero holomorphic function germ and a Weierstrass polynomial of degree $\ell$ in the variable $z_2$.  For the proof of the general case, we consider the germ at $0$ of the holomorphic map $\pi:{\mathbb C}^2\to{\mathbb C}^2$ defined by $(z_1,z_2)\to(\zeta_1,\zeta_2)$.  Since the origin $0$ of ${\mathbb C}^2$ is an isolated point of the common zero-set of $\zeta_1$ and $\zeta_2$, the map $\pi$ is an analytic branched cover (as the germ of a holomorphic map).  Let $C$ be the divisor of $H$ and $\tilde C$ be the image of $C$ under $\pi$ (with multiplicities counted) and let $\tilde H$ be a holomorphic function germ on ${\mathbb C}^2$ at $0$ whose divisor is $\tilde C$.  Since the origin is an isolated point of the common zero-set of $H$ and $\zeta_1$ and since $\zeta_2^\ell$ belongs to the ideal generated by $\zeta_1$ and $H$, it follows that the restriction of $\tilde H$ to $\{\zeta_1=0\}$ is equal to $\zeta_2^{\tilde\ell}$ as holomorphic function germ on ${\mathbb C}$ with $\zeta_2$ as coordinate for some positive integer $\tilde\ell\leq\ell$.  The general case now follows from applying the special case when $H$ is replaced by $\zeta_2^{\tilde\ell-\ell}\tilde H$ and the coordinates $z_1,z_2$ are replaced by $\zeta_1,\zeta_2$.

\bigbreak\noindent(3.3) {\it Step Three.}  Because of Condition (ii) in (3.2), we can now apply the second part of the statement (3.2.2) on Weierstrass polynomials to the case of $H=\hat h_2$ and $\zeta_1=h_1$ and $\zeta_2=w_2$ to get a holomorphic function germ $h_2$ of the form
$$
w_2^{q^2}+\sum_{j=0}^{q^2-1}\theta_j(h_1)w_2^j
$$
which contains $\hat h_2$ as a factor.  (The property of $\hat h_2$ being a reduced holomorphic function germ means that in applying the above statement on Weierstrass polynomial, the divisor $C$ in the proof of the statement (3.2.2) on Weierstrass polynomials is reduced and we do not have to worry about multiplicities of its branches, but this point, though offering some convenience, is not essential.)  Since $h_2$ contains as a factor the multiplier $\hat h_2$ whose order of multiplicity is $\geq\frac{1}{4q}$, it follows $h_2$ is a multiplier whose order of subellipticity is $\geq\frac{1}{4q}$.  Let $h_{2,0}=h_2$ and for $1\leq\nu\leq q^2$ let
$$
h_{2,\nu}=q^2(q^2-1)\cdot(q^2-\nu+1)w_2^{q^2-\nu}+\sum_{j=0}^{q^2-1}j(j-1)\cdots(j-\nu+1)\theta_j(h_1)w_2^{j-\nu},
$$
which is obtained by differentiating $\nu$-times the function $h_2$ with respect to $w_2$ with $h_1$ fixed when $h_2$ is regarded as a function of $h_1$ and $w_2$. Then
$$
dh_{2,\nu}=\eta_\nu dh_1+h_{2,\nu+1}dw_2\quad{\rm for}\ \ 0\leq\nu\leq q^2-1
$$
for a holomorphic function germ $\eta_\nu$ which is the partial derivative of $h_{2,\nu}$ with respect to $h_1$ with $w_2$ fixed when $h_{2,\nu}$ is regarded as a function of $h_1$ and $w_2$.
Let
$$
\tilde h_{2,1}=\left(\frac{\partial h_1}{\partial w_1}\right)_{w_2={\rm const}}
h_{2,1}.$$
From
$$
\begin{aligned}dh_1\wedge dh_2&=dh_1\wedge\left(\eta_0 dh_1+h_{2,1}dw_2\right)\cr
&=dh_1\wedge\,h_{2,1}dw_2\cr
&=\left(\left(\frac{\partial h_1}{\partial w_1}\right)_{w_2={\rm const}}dw_1+
\left(\frac{\partial h_1}{\partial w_2}\right)_{w_1={\rm const}}dw_2\right)\wedge\,h_{2,1} dw_2\cr
&=\tilde h_{2,1} dw_1\wedge dw_2,\cr
\end{aligned}
$$
it follows that $\tilde h_{2,1}$ is a multiplier whose order of subellipticity is $\geq\frac{1}{8q}$.

\medbreak Since $\hat h_2 h_{2,1}$, being a multiple of the multiplier $\hat h_2$, is itself a multiplier, it follows that the linear combination
$$
\alpha \hat h_2 h_{2,1}+\beta\tilde h_{2,1}=\alpha\hat h_2 h_{2,1}+\beta\left(\frac{\partial h_1}{\partial w_1}\right)_{w_2={\rm const}}
h_{2,1}=(h_1)^{3q^2}h_{2,1}
$$
of the two multipliers $\hat h_2 h_{2,1}$ and $\tilde h_{2,1}$ with coefficients, which are holomorphic function germs at $0$ on ${\mathbb C}^2$, is a multiplier whose order of subellipticity is $\geq\frac{1}{8q}$.

\bigbreak\noindent(3.4) {\it Recursive Argument in Step Three.} Now we repeat the above argument with $h_{2,\nu}$ replacing $\hat h_2$ in the following way to conclude by induction on $\nu$ that $(h_1)^{3q^2\nu}h_{2,\nu}$ is a multiplier with order of subellipticity $\geq\frac{1}{2^{\nu+2}q}$ for
$1\leq\nu\leq q^2$.

\medbreak The case of $\nu=1$ was just proved.  Suppose we have proved the step up to some $\nu<q^2$ and we would like to prove the next step of $\nu+1$.
From
$$
\begin{aligned}&dh_1\wedge d\left((h_1)^{3q^2\nu}h_{2,\nu}\right)\cr&=dh_1\wedge
\left(\left((3q^2\nu-1)(h_1)^{3q^2\nu-1}h_{2,\nu}+(h_1)^{3q^2\nu}\eta_\nu\right)dh_1+(h_1)^{3q^2\nu}h_{2,\nu+1}dw_2\right)\cr
&=dh_1\wedge
(h_1)^{3q^2\nu}h_{2,\nu+1}dw_2\cr
&=\left(\left(\frac{\partial h_1}{\partial w_1}\right)_{w_2={\rm const}}dw_1+
\left(\frac{\partial h_1}{\partial w_2}\right)_{w_1={\rm const}}dw_2\right)\wedge
(h_1)^{3q^2\nu}h_{2,\nu+1}dw_2\cr
&=\left(\frac{\partial h_1}{\partial w_1}\right)_{w_2={\rm const}}
(h_1)^{3q^2\nu}h_{2,\nu+1}dw_1\wedge dw_2\cr
\end{aligned}
$$
it follows that $\left(\frac{\partial h_1}{\partial w_1}\right)_{w_2={\rm const}}
(h_1)^{3q^2\nu}h_{2,\nu+1}$ is a multiplier whose order of subellipticity is $\geq\frac{1}{2^{\nu+3}q}$.

\medbreak Since $\hat h_2(h_1)^{3q^2\nu}h_{2,\nu+1}$, being a multiple of the multiplier $\hat h_2$, is itself a multiplier, it follows that the linear combination
$$
\begin{aligned}&\alpha \hat h_2(h_1)^{3q^2\nu}h_{2,\nu+1}+\beta\left(\frac{\partial h_1}{\partial w_1}\right)_{w_2={\rm const}}
(h_1)^{3q^2\nu}h_{2,\nu+1}\cr
&=\left(\alpha\hat h_2+\beta\left(\frac{\partial h_1}{\partial w_1}\right)_{w_2={\rm const}}\right)
(h_1)^{3q^2\nu}h_{2,\nu+1}\cr
&=(h_1)^{3q^2(\nu+1)}h_{2,\nu+1}\cr
\end{aligned}
$$
of the two multipliers $\hat h_2(h_1)^{3q^2\nu}h_{2,\nu+1}$ and $\beta\left(\frac{\partial h_1}{\partial w_1}\right)_{w_2={\rm const}}
(h_1)^{3q^2\nu}h_{2,\nu+1}$ with coefficients, which are holomorphic function germs at $0$ on ${\mathbb C}^2$, is a multiplier with order of subellipticity $\geq\frac{1}{2^{\nu+3}q}$.  This finishes the proof by induction on $\nu$ that $(h_1)^{3q^2\nu}h_{2,\nu}$ is a multiplier with order of subellipticity $\geq\frac{1}{2^{\nu+2}q}$ for
$1\leq\nu\leq q^2$.  Since $h_{2,q^2}$ is equal to $(q^2)!$, it follows that $(h_1)^{3q^4}$ is a multiplier whose order of subellipticity is $\geq\frac{1}{2^{q^2+2}q}$.
By the real radical property of multipliers we conclude that $h_1$ is a multiplier whose order of subellipticity is $\geq\frac{1}{3q^5 2^{q^2+2}}$.

\medbreak Since by Condition (ii) of (3.2) the ideal generated by $h_1$ and $\hat h_2$ contains the $q^2$-th power of the maximum ideal ${\mathfrak m}_{{\mathbb C}^2,0}$ of ${\mathbb C}^2$ at $0$, it follows that $w_1$ and $w_2$ are multipliers whose order of subellipticity is $\geq\frac{1}{3q^7 2^{q^2+2}}$. By taking the Jacobian determinant of $w_1$ and $w_2$, we end up with $1$ being a multiplier whose order of subellipticity is $\frac{1}{3q^7 2^{q^2+3}}$.

\bigbreak\noindent(3.5) {\it Remark.}  For use in (4.3) below we would like to remark that the above arguments work in the same way when the holomorphic function germ $h_2$ is chosen to be
$$
(h_1)^r\left(w_2^{q^2}+\sum_{j=0}^{q^2-1}\theta_j(h_1)w_2^j\right)
$$
instead of
$$
w_2^{q^2}+\sum_{j=0}^{q^2-1}\theta_j(h_1)w_2^j
$$
if $r$ is effective in the sense that $r$ is bounded by an explicit function of $q$. Of course, the effective lower bound of the order of subellipticity of the constant function $1$ as a multiplier needs to be correspondingly modified to be $\frac{1}{3q^7 2^{q^2+r+3}}$.

\bigbreak\noindent{\bf\S 4.} {\sc Effective Kohn Algorithm Applied to Catlin-D'Angelo's Example}

\bigbreak We now apply the algebraic geometric techniques in the effective Kohn algorithm to the example of Catlin-D'Angelo given in [Catlin-D'Angelo2010] for which the full-real-radical Kohn algorithm is ineffective.

\bigbreak\noindent(4.1) {\it Catlin-D'Angelo's Example of Ineffectiveness of Full-Real-Radical Kohn Algorithm.}  Let $K >M\geq 2$ and $N\geq 3$.  The special domain $\Omega$
$$
{\rm Re}\,w+|F_1(z_1,z_2)|^2+|F_2(z_1,z_2)|^2<0
$$
in ${\mathbb C}^3$ is defined by the two holomorphic functions $F_1(z_1,z_2)=z_1^M$ and $F_2=z_2^N+z_2z_1^K$ on ${\mathbb C}^2$.  The origin of ${\mathbb C}^3$ is the boundary point of $\Omega$ whose scalar and vector multipliers we consider.
The following is reproduced from pp.81-82 of [Catlin-D'Angelo2010] in the notations and terminology used in this note.  By Weierstrass division (applied to the Weierstrass polynomial which is the product $F_2$ and a nowhere holomorphic function germ), modulo $F_2$ every element of ${\mathcal O}_{{\mathbb C}^2,0}$ is equal to
$$
\sum_{j=0}^{N-1}a_j(z_1)z_2^j
$$
for some holomorphic function germs $a_0,\cdots,a_{N-1}$ on ${\mathbb C}^2$ at $0$.  Modulo $F_1$ we can write
$$
a_j(z_1)=\sum_{k=0}^{M-1}b_{jk}z_1^k
$$
for $0\leq j\leq N-1$, where $b_{jk}\in{\mathbb C}$ for $0\leq j\leq N-1,\,0\leq k\leq M-1$.  Hence the multiplicity $q$ given by
$$
\dim_{\mathbb C}\left({\mathcal O}_{{\mathbb C}^2,0}\left/\sum_{j=1}^2
{\mathcal O}_{{\mathbb C}^2,0} F_j\right.\right)
$$
is $\leq MN$.

\medbreak The full-real-radical Kohn algorithm proceeds as follows in this example.  Let $g$ be the Jacobian determinant
$$\frac{\partial(F_1,F_2)}{\partial(z_1,z_2)}=\det\left(\begin{matrix}z_1^{M-1}&0\cr
Kz_2z_1^{K-1}&Nz_2^{N-1}+z_1^K\cr
\end{matrix}\right)=Nz_1^{M-1}z_2^{N-1}+z_1^{K+M-1}$$
of $F_1,F_2$ with respect to $z_1,z_2$.  We use the notations in (2.10).  The ideal $J_0$ is the principal ideal with the irreducible function germ $g$ as the generator and its radical $I_0$ is the same as $J_0$. The ideal $J_1$ is
$${\mathcal O}_{{\mathbb C}^2,0}g+\sum_{j=1}^2{\mathcal O}_{{\mathbb C}^2,0}\frac{\partial(F_j,g)}{\partial(z_1,z_2)},$$
where
$$\frac{\partial(F_1,g)}{\partial(z_1,z_2)}=\det\left(\begin{matrix}z_1^{M-1}&0\cr
\frac{\partial g}{\partial z_1}&N(N-1)z_1^{M-1}z_2^{N-2}\cr
\end{matrix}\right)=N(N-1)z_1^{2M-2}z_2^{N-2}.$$
Since $g^2$ modulo $\frac{\partial(F_1,g)}{\partial(z_1,z_2)}$ is equal to $z_1^{2(K+M-1)}$, it follows from $K+M-1\geq 1$ that the holomorphic function germ $z_1$ at $0$ belongs to the radical $I_1$ of $J_1$.  From $N\geq 3$ we conclude that modulo $(z_2)^2$ the three holomorphic function germs $g, F_1, F_2$ become respectively $z_1^{K+M-1}$, $z_1^M$, $z_2z_1^K$.  Hence the ideal $J_1$ generated by the three Jacobian determinants formed from pairs out of $g, F_1, F_2$ is contained in the ideal generated by $z_1^{M+K-2}$ and $z_2$.  This means that $z_1^m$ cannot be in $J_1$ for $m<M+K-2$, otherwise $z_1^m$ belongs to the ideal generated by $z_1^{M+K-2}$ and $z_2$, which is a contradiction.  Since the holomorphic function germ $z_1$ belongs to $I_1$, this means that the smallest positive integer $p_1$ satisfying $(I_1)^{p_1}\subset J_1$ must be $\geq M+K-2\geq K$.  Thus, the full-real-radical Kohn algorithm is not effective, because $K$ is arbitrary and there cannot be any function of $NM$ which bounds $K$.

\bigbreak\noindent(4.2) {\it Remark.} In (4.1) when we carry out the full-real-radical Kohn algorithm for Catlin-D'Angelo's example, we stopped after showing the algorithm to be ineffective.  For later comparison, we now carry out the remaining steps of the algorithm until we produce the constant function $1$ as a multiplier.  We have seen that $z_1$ belongs to $I_1$.  Since all three holomorphic function germs $g$, $\frac{\partial g}{\partial z_1}$ and $\frac{\partial g}{\partial z_2}$ contain $z_1$ as a factor, it follows that $J_1$ is contained in the principal ideal generated by $z_1$ and $I_1$ must be equal to the principal ideal generated by $z_1$.  The function germ $z_2^N=F_2-z_2z_1^K$ is a pre-multiplier in $V_2$.  The Jacobian determinant $\frac{\partial(z_1,z_2^N)}{\partial(z_1,z_2)}=Nz_2^{N-1}$ belongs to $J_2$.  Hence $z_2$ belongs $I_2$ and we can conclude that $I_2$ is the maximum ideal ${\mathfrak m}_{{\mathbb C}^2,0}$ of ${\mathbb C}^2$ at $0$.  By taking the Jacobian determinant of the elements $z_1, z_2$ of $I_2$, we conclude that $1$ is a multiplier.  To get to the multiplier $1$ from $F_1,F_2$, we have to perform differentiation $4$ times in the construction of Jacobian determinants.

\bigbreak\noindent(4.3) {\it Effective Kohn Algorithm for Catlin-D'Angelo's Example.}
We now carry out concretely the steps in the effective Kohn algorithm for Catlin-D'Angelo's example to illustrate the difference between the full-real-radical Kohn algorithm and
the effective Kohn algorithm.

\medbreak The key point in the effective Kohn algorithm is to construct a Weierstrass polynomial $h_2$ in one coordinate $w_2$ such that (i) $h_2$ contains as a factor a multiplier $\hat h_2$ which is obtained in a procedure involving the Jacobian determinant of two ${\mathbb C}$-linear combinations of the defining holomorphic functions $F_1,\cdots,F_N$ of the special domain and (ii) the coefficients of $h_2$ are holomorphic function germs of some ${\mathbb C}$-linear combination $h_1$ of $F_1,\cdots,F_N$.  Then by using induction on $\nu$ we show, with effectiveness, that the Jacobian determinant of $h_1$ and the function $(h_1)^{m_\nu}\left(\frac{\partial^\nu h_2}{\partial w_2^\nu}\right)_{h_1={\rm constant}}$ is a multiplier (for some effective positive integer $m_\nu$), resulting finally in the conclusion that $h_1$ is a multiplier.  In the key argument the Weierstrass polynomial $h_2$ can be replaced by the product of an effective power of $h_1$ and a Weierstrass polynomial (see Remark (3.4)).

\medbreak In the example of Catlin-D'Angelo where the defining functions for the special domain in ${\mathbb C}^3$ are the two holomorphic functions $F_1(z_1,z_2)=z_1^M$ and $F_2=z_2^N+z_2z_1^K$ on ${\mathbb C}^2$, we can use $h_1=F_1$ and use
$$\frac{z_1}{N}\frac{\partial(F_1,F_2)}{\partial(z_1,z_2)}=z_1^M z_2^{N-1}+\frac{1}{N}\,z_1^{K+M}$$
as $h_2$ which is the product of $h_1$ and the Weierstrass polynomial
$$
z_2^{N-1}+\frac{1}{N}\,z_1^K=z_2^{N-1}+\frac{1}{N}(h_1)^m
$$
in the variable $z_2$, where $m=\frac{K}{M}$ which we assume for the time being to be a positive integer. It turns out that the argument used in the effective Kohn algorithm works in the same way without the assumption that $m=\frac{K}{M}$ is a positive integer.  This assumption used in the setup merely motivates the steps of the argument.

\medbreak For Catlin-D'Angelo's example, the induction on $\nu$ to show, with effectiveness, that the Jacobian determinant of $h_1$ and the function $(h_1)^{m_\nu}\left(\frac{\partial^\nu h_2}{\partial w_2^\nu}\right)_{h_1={\rm constant}}$ is a multiplier (for some effective positive integer $m_\nu$) is translated (after obvious modifications) to verifying by induction on $j$ that each $H_j$ defined by $H_j=z_1^{(j+1)(M-1)}z_2^{N-j}$ is a multiplier for $1\leq j\leq N$, because the $\nu$-th derivative of
$$
z_2^{N-1}+\frac{1}{N}\,z_1^K=z_2^{N-1}+\frac{1}{N}(h_1)^m
$$
with respect to $z_2$ with $h_1$ being kept constant is $(N-1)\cdots(N-\nu)z_2^{N-1-\nu}$ for $1\leq\nu\leq N-1$.

\medbreak At this point we can forget that the use of $H_j=z_1^{(j+1)(M-1)}z_2^{N-j}$ for $1\leq j\leq N$ is motivated by the steps of the effective Kohn algorithm.  We now simply carry out the induction on $j$ to verify that $H_j=z_1^{(j+1)(M-1)}z_2^{N-j}$ is a multiplier whose order of subellipticity is $\geq\frac{1}{2^{j+2}}$ for $1\leq j\leq N$.

\medbreak Since the Jacobian determinant
$$
g=\frac{\partial(F_1,F_2)}{\partial(z_1,z_2)}=Nz_1^{M-1}z_2^{N-1}+z_1^{K+M-1}
$$
is multiplier whose order of subellipticity is $\geq\frac{1}{4}$, the Jacobian determinant
$$
\frac{\partial(F_1,g)}{\partial(z_1,z_2)}=N(N-1)z_1^{2M-2}z_2^{N-1}
$$
is a multiplier whose order of subellipticity is $\geq\frac{1}{8}$, which means that $H_1$ is a multiplier whose order of subellipticity is $\geq\frac{1}{8}$.  Suppose $H_j$ has been verified to be a multiplier whose order of subellipticity is $\geq\frac{1}{2^{j+2}}$ for some $1\leq j<N$.  Then the Jacobian determinant
$$
\begin{aligned}\frac{\partial(F_1,H_j)}{\partial(z_1,z_2)}&=
\det\left(\begin{matrix}Mz_1^{M-1}&0\cr \frac{\partial H_j}{\partial z_1}&(N-j)z_1^{(j+1)(M-1)}z_2^{N-j-1}\cr\end{matrix}\right)\cr
&=M(N-j)z_1^{(j+2)(M-1)}z_2^{N-j-1}\cr
\end{aligned}
$$
is a multiplier whose order of subellipticity is $\geq\frac{1}{2^{j+3}}$, which means that $H_{j+1}$ is a multiplier whose order of subellipticity is $\geq\frac{1}{2^{j+3}}$.  This finishes the induction argument and we know that $H_N=z_1^{(N+1)(M-1)}$ is a multiplier whose order of subellipticity is $\geq\frac{1}{2^{N+2}}$.

\medbreak By taking the effective $(N+1)(M-1)$-th root, we conclude that the holomorphic function germ $z_1$ is a multiplier whose order of subellipticity is $\geq\frac{1}{2^{N+2}(N+1)(M-1)}$ and the holomorphic function germ $z_2z_1^K$ which contains $z_1$ as a factor is a multiplier whose order of subellipticity is $\geq\frac{1}{2^{N+2}(N+1)(M-1)}$.  Then
$$z_2^N=(z_2^N+z_2z_1^K)-z_2z_1^K=F_2-z_2z_1^K$$ is a pre-multiplier whose differential has order of subellipticity $\geq\frac{1}{2^{N+3}(N+1)(M-1)}$.   Since both $z_1$ and $z_2^N$ are pre-multipliers whose differentials have order of subellipticity $\geq\frac{1}{2^{N+3}(N+1)(M-1)}$, the Jacobian determinant
$$\frac{\partial(z_1,z_2^N)}{\partial(z_1,z_2)}=\det\left(\begin{matrix}1&0\cr 0&Nz_2^{N-1}\cr\end{matrix}\right)=Nz_2^{N-1}
$$
is a multiplier whose order of subellipticity is $\geq\frac{1}{2^{N+3}(N+1)(M-1)}$, which means that $z_2^{N-1}$ is a multiplier whose order of subellipticity is $\geq\frac{1}{2^{N+3}(N+1)(M-1)}$.
By taking the $(N-1)$-th root of $z_2^{N-1}$, we conclude that $z_2$ is a multiplier whose order of subellipticity is $\geq\frac{1}{2^{N+3}(N+1)(N-1)(M-1)}$.  Finally $1$ which is the Jacobian determinant $\frac{\partial(z_1,z_2)}{\partial(z_1,z_2)}$ is a multiplier whose order of subellipticity is $\geq\frac{1}{2^{N+4}(N+1)(N-1)(M-1)}$.  This algorithm is effective, because $\frac{1}{2^{N+4}(N+1)(N-1)(M-1)}$ is bounded from below by the explicit function $\frac{1}{2^{q+4}(q+1)(q-1)^2}$ of $q=NM$.

\bigbreak\noindent(4.4) {\it Remark.} In carrying out above the effective Kohn algorithm for Catlin-D'Angelo's example,
to get to the multiplier $1$ from $F_1,F_2$, we have to perform differentiation $N+2$ times in the construction of Jacobian determinants.  Compared to the full-real-radical Kohn algorithm which
requires only $4$ differentiation to terminate, to avoid ineffectiveness in the taking of roots in the effective Kohn algorithm we choose the option of performing more, but still an effective number of, differentiations.

\bigbreak\noindent(4.5) {\it Geometric Reason for Ineffectiveness of Full-Real-Radical Kohn Algorithm for Catlin-D'Angelo's Example.} The above discussion shows by computation why in Catlin-D'Angelo's example the full-real-radical Kohn algorithm is ineffective while the effective Kohn algorithm gives effectiveness.  Now we would like to analyze geometrically why such a phenomenon occurs.
When $F_2=z_2^N+z_2z_1^K$ is regarded as a polynomial in $z_2$, its degree $N$ is effective (in the sense of being bounded by an explicit function of $q=MN$) but its discriminant obtained by eliminating $z_2$ from $F_2$ and $\frac{\partial F_2}{\partial z_2}$, as a function germ in $z_1$ vanishes to an order at  $z_1=0$ which is a function of $K$ and is not effective.  In other words, the $N$ roots (in $z_2$) of $F_2=0$ as $N$ functions of $z_1$ are close together near $z_1=0$ to an order which is a function of $K$ and is not effective.  The discriminant of $F_2$ and the closeness of the $N$ roots of $F_2=0$ enter the picture, because $F_1=z_1^M$ depends only on $z_1$ and the Jacobian determinant of $F_1, F_2$ is the first multiplier in the algorithm.  Because of the ineffectiveness of the vanishing order of the discriminant of $F_2$ at a function in $z_1$ at $z_1=0$, the step of root-taking is ineffective.  On the other hand, the effective Kohn algorithm replaces ineffective root-taking of the discriminant of a Weierstrass polynomial by differentiating the Weierstrass polynomial with respect to its variable as many times as its degree to avoid the ineffective root-taking.

\bigbreak\noindent{\bf\S 5.} {\sc Multipliers in More General Setting}

\bigbreak We now discuss the generalization of Kohn's technique of multipliers to more general systems of partial differential equations.

\bigbreak\noindent(5.1) {\it Generalization of Kohn's Technique of Multiplier Ideal Sheaves to More General Setting.}  The generalization of Kohn's technique of multiplier ideal sheaves to a more general setting comes from looking at Kohn's technique for the complex Neumann problem from the following perspective.  The subelliptic estimate for a bounded smooth weakly pseudoconvex domain $\Omega$ in ${\mathbb C}^n$ at its boundary point $P$ seeks to estimate
$$\left\|\left|\varphi\right|\right\|_\varepsilon^2=\left\|\Lambda^\varepsilon\varphi\right\|^2$$
by a constant $C_\varepsilon$ times
$$
Q(\varphi,\varphi)=\left\|\bar\partial\varphi\right\|^2+\left\|\bar\partial^*\varphi\right\|^2+\left\|\varphi\right\|^2,
$$
for some $\varepsilon>0$, for all smooth test $(0,1)$-form $\varphi$ on $U\cap\bar\Omega$ with compact support which belongs to the domain of the actual adjoint $\bar\partial^*$ of $\bar\partial$ (where $U$ is an open neighborhood of $P$ in ${\mathbb C}^n$).  For the convenience of discussion, we simply say that $\Lambda^\varepsilon\varphi$ is {\it estimable} on $U$ when
$$
\left\|\Lambda^\varepsilon\varphi\right\|^2\leq C_\varepsilon Q(\varphi,\varphi).
$$
In general, we say that some expression $\psi$ defined from $\varphi$ is  {\it estimable} on $U$ (or simply {\it estimable}) if
$$
\left\|\psi\right\|^2\leq CQ(\varphi,\varphi)
$$
for some constant $C$ independent of $\varphi$ (which is smooth on $U\cap\partial\Omega$ with compact support).
The starting point is the basic identity
$$
\left\|\bar\partial\varphi\right\|^2+\left\|\bar\partial^*\varphi\right\|^2=\left\|\bar\nabla\varphi\right\|^2+\int_{\partial\Omega}\left<{\rm Levi}_{\partial\Omega},\,\bar\varphi\wedge\varepsilon\right>,
$$
where $\bar\nabla$ is the (covariant) differentiation of $\varphi$ in the $(0,1)$-direction and ${\rm Levi}_{\partial\Omega}$ is the Levi form $\partial\bar\partial r$ of the boundary $\partial\Omega$ of $\Omega$ when $\Omega$ is locally defined by $r<0$ with $dr\equiv1$ on $\partial\Omega$.  In particular,
$$
\left\|\bar\nabla\varphi\right\|^2\leq Q(\varphi,\varphi)\quad{\rm and}\quad
\left\|\bar\partial^*\varphi\right\|^2\leq Q(\varphi,\varphi).
$$
Together with $\left\|\varphi\right\|^2\leq Q(x,y)$, this means that both $\bar\nabla\varphi$ and $\bar\partial^*\varphi$, as well as $\varphi$, are estimable.  The expressions
$$
\bar\nabla\varphi=\left(\bar\partial_j\varphi_{\bar k}\right)_{1\leq j,k\leq n}\quad{\rm and}\quad\bar\partial^*\varphi=-\sum_{j=1}^n\partial_j\varphi_{\bar j}
$$
are linear combinations of first-order partial derivatives of the components of $\varphi$.  A multiplier $F$ means the estimability of $\Lambda^\varepsilon(F\varphi)$ and a vector multiplier $\theta$ means the estimability of $\Lambda^\varepsilon(\theta\cdot\varphi)$.  Kohn's technique is to use the estimability of $\bar\nabla\varphi$, $\bar\partial^*\varphi$ and $\varphi$ and apply algebraic manipulations and integration by parts to construct from the estimability of $\Lambda^\varepsilon(F\varphi)$ and $\Lambda^\varepsilon(\theta\cdot\varphi)$ other $F^\prime$ and $\theta^\prime$ with estimable $\Lambda^\varepsilon(F^\prime\varphi)$ and $\Lambda^\varepsilon(\theta^\prime\cdot\varphi)$.  For such manipulations it does not matter what the meaning of $Q(\varphi,\varphi)$ is.  Moreover, the operations of integration by parts are along the tangent directions of the boundary, because $\Lambda^\varepsilon$ is the pseudo-differential operator corresponding to the $\left(\frac{\varepsilon}{2}\right)$-th power of $1$ plus the Laplacian in tangential coordinates of the boundary.

\bigbreak\noindent(5.1.1) For our generalization of the technique of multiplier ideal sheaves, we use the following simple setting which highlights the core argument of the technique.  Fix an integer $q\geq 2$.  Let $\Omega$ be an open neighborhood of $0$ in ${\mathbb R}^m$ and $Y_{j\,\nu}$ be complex-valued smooth differential operators on $\Omega$ for $1\leq j\leq N$ and $1\leq\nu\leq q$.  For any $q$-tuple $\varphi=\left(\varphi_1,\cdots,\varphi_q\right)$ of smooth complex-valued functions with compact support on $\Omega$, let
$$
Q(\varphi,\varphi)=\left\|\varphi\right\|^2+\sum_{1\leq j\leq N,\,1\leq\nu\leq q}\left\|Y_{j\,\nu}\varphi_\nu\right\|^2,
$$
where $\left\|\cdot\right\|$ means the $L^2$ norm on $\Omega$.  An expression $\psi$ of $\varphi$ of the form $\sum_{\nu=1}^q Z_\nu\varphi_\nu$ (where each $Z_\nu$ is a pseudo-differential operator on $\Omega$) is said to be {\it estimable} on an open neighborhood $U$ of $0$ in $\Omega$ (or simply {\it estimable}) if there is a positive constant $C$ such that
$$
\left\|\psi\right\|^2\leq CQ(\varphi,\varphi)
$$
for all $q$-tuple $\varphi=\left(\varphi_1,\cdots,\varphi_q\right)$ of smooth complex-valued test functions with compact support on $U$.  When $\psi$ is vector-valued instead of scalar-valued, the estimability of $\psi$ on $U$ means the estimability of each of its components on $U$.  When we have two such expressions $\psi$ and $\tilde\psi$, we say that the inner product $\left(\psi,\tilde\psi\right)$ is {\it estimable} on $U$ if
$$
\left(\psi,\,\tilde\psi\right)\leq CQ(\varphi,\varphi)
$$
all $q$-tuple $\varphi=\left(\varphi_1,\cdots,\varphi_q\right)$ of smooth complex-valued test functions with compact support on $U$.
We refer to $C$ as the constant of estimability of $\psi$ or $(\psi,\tilde\psi)$.

\bigbreak\noindent(5.1.2) Let $\Lambda^\varepsilon$ be the the pseudo-differential operator which is the $\left(\frac{\varepsilon}{2}\right)$-th power of $1$ plus the Laplacian in the coordinates of ${\mathbb R}^m$. We introduce three kinds of multipliers: (i) scalar multiplier, (ii) vector multiplier, and (iii) matrix multiplier.

\medbreak The germ at $0$ of a smooth function $\alpha$ is a {\it scalar multiplier} at $0$ with {\it order of subellpticity} $\geq\varepsilon$ if $\Lambda^\varepsilon(\alpha\varphi_\nu)$ (for $1\leq\nu\leq q$) is estimable on $U$ for some $\varepsilon>0$ and some open neighborhood $U$ of $0$ in $\Omega$ on which $\alpha$ is defined.  Clearly the product of a scalar multiplier with any smooth function is again a scalar multiplier with no change in the order of subellipticity.  By considering the commutator $\left[\Lambda^\varepsilon,\,\alpha\right]$ of pseudo-differential operators, we conclude that for $0<\varepsilon\leq 1$ the estimatbility of $\Lambda^\varepsilon(\alpha\varphi_\nu)$ on $U$ is equivalent to the estimability of $\alpha(\Lambda^\varepsilon\varphi_\nu)$ on $U$, because $\left\|\varphi\right\|^2$ is estimable on $U$.

\medbreak The germ at $0$ of a smooth $q$-tuple of smooth complex-valued functions $\vec a=(a_1,\cdots,a_q)$ is called a {\it vector multiplier} at $0$ with {\it order of subellpticity} $\geq\varepsilon$ if $\Lambda^\varepsilon(\sum_{\nu=1}^q a_\nu\varphi_\nu)$ is estimable on $U$ for some $\varepsilon>0$ and some open neighborhood $U$ of $0$ in $\Omega$ on which $\vec a$ is defined.  Clearly the product of a vector multiplier with any smooth function is again a vector multiplier with no change in the order of subellipticity. Again, for $0<\varepsilon\leq 1$ the estimatbility of $\Lambda^\varepsilon(\sum_{\nu=1}^q a_\nu\varphi_\nu)$ on $U$ is equivalent to the estimability of $\sum_{\nu=1}^q a_\nu(\Lambda^\varepsilon\varphi_\nu)$ on $U$.

\medbreak   An $q\times q$ matrix ${\mathbf a}=\left(a_{jk}\right)_{1\leq j,\ell\leq q}$ is called a {\it matrix multiplier} at $0$ with {\it order of subellpticity} $\geq\varepsilon$ if every one of its rows $\vec a_j=(a_{j1},\cdots,a_{jq})$ (for $1\leq r\leq q$) is a vector multiplier at $0$ with order of subellipticity $\geq\varepsilon$.  Clearly a matrix multiplier multiplied on the left by a $q\times q$ matrix with smooth functions as entries yields a matrix multiplier with no change in the order of subellipticity.

\medbreak Some simple relations among scalar multipliers, vector multipliers, and matrix multipliers are as follows.  The product of a scalar multiplier with any row $q$-vector with smooth functions as components is a vector multiplier.  The product of a scalar multiplier with any $q\times q$ matrix with smooth functions as entries is a matrix multiplier.  Any vector multiplier (as a row vector) multiplied on the left by a column $q$-vector with smooth functions as components yields a matrix multiplier.  By Cramer's rule the determinant of a matrix multiplier is a scalar multiplier. Any matrix multiplier multiplied on the left by a row $q$-vector with smooth functions as components yields a vector multiplier.

\bigbreak\noindent(5.1.3) Just like the real radical property (C) of Kohn's multipliers in (2.1), scalar multipliers here enjoy the same real radical property that if $\alpha$ is a scalar multiplier at $0$ with order of subellipticity $\geq\varepsilon$ (for some $0<\varepsilon\leq 1$) and $\beta$ is a smooth complex-valued function germ at $0$ such that $\left|\beta\right|^\sigma\leq\left|\alpha\right|$ for some $\sigma\in{\mathbb N}$, then $\beta$ is a scalar multiplier with order of subellipticity $\geq\frac{\varepsilon}{\sigma}$.  The proof is completely analogous to the proof of [Kohn1979, p.98, Lemma 4.3.4] and is as follows.  Let $U$ be an open neighborhood of $0$ in $\Omega$ such that both $\alpha$ and $\beta$ are represented by smooth functions on $U$ and $\varphi$ is a test $q$-tuple of smooth functions on $U$ with compact support.  Let $\eta=\frac{\varepsilon}{\sigma}$.  Since
$$
\begin{aligned}\left\|\Lambda^{\sigma\eta}(\beta^\sigma\varphi)\right\|_{L^2(U)}&\leq C_1^*\left\|\beta^\sigma\Lambda^{\sigma\eta}\varphi\right\|_{L^2(U)}+C_1^{**}\left\|\varphi\right\|_{L^2(U)}\cr
&\leq C_2^*\left\|\alpha\Lambda^{\sigma\eta}\varphi\right\|_{L^2(U)}+C_2^{**}\left\|\varphi\right\|_{L^2(U)}\cr
&\leq C_3^*\left\|\Lambda^{\sigma\eta}(\alpha)\varphi\right\|_{L^2(U)}+C_3^{**}\left\|\varphi\right\|_{L^2(U)},\cr
\end{aligned}
$$
it suffices to prove the statement that
$$
\left\|\Lambda^{\tau\eta}(\beta^\tau\varphi)\right\|_{L^2(U)}
\leq C_\tau\left\|\Lambda^{\sigma\eta}(\beta^{\sigma\eta}\varphi)\right\|_{L^2(U)}+C_\tau^\prime\left\|\varphi\right\|_{L^2(U)}\leqno{(5.1.3.1)_\tau}
$$
for $1\leq\tau\leq\sigma$.
The statement $(5.1.3.1)_\tau$ follows from descending induction on $\tau$ for $1\leq\tau\leq\sigma$, because
$$
\begin{aligned}\left\|\Lambda^{\tau\eta}(\beta^\tau\varphi)\right\|_{L^2(U)}^2
&=
\left(\Lambda^{(\tau+1)\eta}(\beta^{\tau+1}\varphi),\Lambda^{(\tau-1)\eta}(\beta^{\tau-1}\varphi))\right)_{L^2(U)}
+\hat C\left\|\varphi\right\|^2\cr
&\leq\left\|\Lambda^{(\tau+1)\eta}(\beta^{\tau+1}\varphi)\right\|_{L^2(U)}
\left\|\Lambda^{(\tau-1)\eta}(\beta^{\tau-1}\varphi)\right\|_{L^2(U)}
+\tilde C\left\|\varphi\right\|^2\cr
\end{aligned}
$$
for $1\leq\tau<\sigma$, where $C_j^*, C_j^{**}, C_\tau, C_\tau^\prime, \hat C, \tilde C$ are constants independent of $\varphi$.

\medbreak In particular, if $\alpha$ is a scalar multiplier at $0$ with order of subellipticity $\geq\varepsilon$ (for some $0<\varepsilon\leq 1$), then its complex-conjugate $\bar\alpha$ is also scalar multiplier at $0$ with order of subellipticity $\geq\varepsilon$.

\bigbreak A very important part of the technique of multiply ideal sheaves is the differential relations among the scalar multipliers, vector multipliers, and matrix multipliers.  This is represented by two procedures involving differentiation.  The first procedure produces a new vector multiplier from a matrix multiplier.  The second procedure produces a new vector multiplier from a scalar multiplier.  The second procedure is similar to the procedure (B)(i) for Kohn's multipliers for the complex Neumann problem in (2.1).  The first procedure is a new one, even in the special case of Kohn's multipliers for the complex Neumann problem.  The following theorem presents a unified version of both procedures, which yields both procedures as special cases.

\bigbreak\noindent(5.2) {\it Theorem (Generation of Vector Multiplier from Matrix Multiplier or Scalar Multiplier).}  Let $X_1,\cdots,X_q$ be complex-valued smooth first-order differential operators on $\Omega$ whose adjoint operators are $X_1^*,\cdots,X_q^*$ with respect to the $L^2$ inner product on $\Omega$ such that each $X_j^*\varphi$ is estimable on $\Omega$ for $1\leq j\leq q$.  Let $\Gamma_{k\ell}$ be a smooth complex-valued function on $\Omega$ for $1\leq k,\ell\leq q$ such that $\sum_{1\leq k,\ell\leq q}\Gamma_{k\ell}X_k\varphi_\ell$ is estimable on $\Omega$.  Let $\varepsilon_1,\,\varepsilon_2$ be positive numbers $\leq 1$.  Let ${\mathbf a}=\left(a_{jk}\right)_{1\leq j,k\leq q}$ be a matrix of multipliers at $0$ so that each of its rows $\vec a_j=\left(a_{jk}\right)_{1\leq k\leq q}$ is a vector multiplier at $0$ with order of subellipticity $\geq\varepsilon_1$ for $1\leq j\leq q$.  Let $\alpha$ be a scalar multiplier at $0$ with order of subellipticity $\geq\varepsilon_2$.  Let $\left(A_{jk}\right)_{1\leq j,k\leq q}$ be a matrix of smooth complex-valued function germs at $0$ such that $\sum_{\ell=1}^q A_{j\ell}a_{\ell k}$ equals to the Kronecker delta $\delta_{jk}$ times $\alpha$ for $1\leq j,k\leq q$.  Let
$$
b_j=\sum_{1\leq p,\ell,k\leq q}\Gamma_{pk}A_{k\ell}(X_pa_{\ell\,j})
$$
and $\vec b=\left(b_j\right)_{1\leq j\leq q}$.  Then $\vec b$ is a vector multiplier at $0$ whose order of subellipticity is $\geq\frac{1}{2}\min(\varepsilon_1,\varepsilon_2)$.  In particular, the following two special cases hold.

\medbreak\noindent(i) For $1\leq j\leq q$ let
$$c_j=\sum_{1\leq p,\ell\leq q}{\rm adj}({\mathbf a})_{p\ell}(X_p a_{\ell\,j})$$
and $\vec c=\left(c_j\right)_{1\leq j\leq q}$, where ${\rm adj}({\mathbf a})$ is the adjoint matrix of ${\mathbf a}$ (so that the matrix product of ${\rm adj}({\mathbf a})$ and ${\mathbf a}$ is equal to $\det({\mathbf a})$ times the identity matrix of order $q$).  Then $\vec c$ is a vector multiplier at $0$ whose order of subellipticity is $\geq\frac{\varepsilon_1}{2}$.

\medbreak\noindent(ii) Let $d_j=\sum_{k=1}^q\Gamma_{kj}X_k\alpha$ for $1\leq j\leq q$ and $\vec d=\left(d_j\right)_{1\leq j\leq q}$.  Then $\vec d$ is a vector multiplier at $0$ whose order of subellipticity is $\geq\frac{\varepsilon_2}{2}$.

\bigbreak\noindent(5.2.1) {\it Proof of Theorem (5.2).} Let $\varepsilon=\min(\varepsilon_1,\varepsilon_2)$.  Let $U$ be an open neighborhood of $0$ in $\Omega$ such that $\alpha$ and the vector multiplier $\vec a_j=\left(a_{jk}\right)_{1\leq k\leq q}$ (for $1\leq j\leq q$) are defined and smooth on $U$ and
$\Lambda^\varepsilon(\sum_{\nu=1}^q a_{j\nu}\varphi_\nu)$ is estimable on $U$ for smooth test functions $\varphi=(\varphi_1,\cdots,\varphi_q)$ on $U$ with compact support.
Let $\psi$ be a linear combination of $\varphi_j$ (for $1\leq j\leq q$) with smooth functions on $U$ as coefficients and which we will specify more precisely later.  By the Cauchy-Schwarz inequaltiy, for $1\leq p,\ell\leq q$ the inner product
$$
\left(\Lambda^\varepsilon\left(\sum_{j=1}^q a_{\ell\,j}\,\varphi_j\right),\,X_p^*\psi\right)
$$
is estimable on $U$, because $\vec a_\ell$ is a vector-multiplier for $1\leq\ell\leq q$ and $X_p^*\psi$ is estimable (from the estimability of $X_p^*\varphi$) for $1\leq p\leq q$.   Note that the constant of estimability depends on the smooth coefficient functions in the linear combination $\psi$ of $\varphi_1,\cdots,\varphi_q$ which are yet to be specified.  Integration by parts applied to $X_p$ (by switching $X_p$ over to $X_p^*$ in the inner product) yields the estimability of
$$
\left(\Lambda^\varepsilon\left(\sum_{j=1}^q (X_p a_{\ell\,j})\varphi_j\right),\,\psi\right)
+\left(\Lambda^\varepsilon\left(\sum_{j=1}^q a_{\ell\,j}(X_p \varphi_j)\right),\,\psi\right)\leqno{(5.2.1.1)}
$$
on $U$ for $1\leq p,\ell\leq q$ after we take care of the error terms from the commutator of pseudodifferential operators in the standard way.

\medbreak Now we apply $\sum_{\ell=1}^q A_{k\ell}$ to (5.2.1.1) to get the estimability of
$$
\left(\Lambda^\varepsilon\left(\sum_{\ell,j=1}^qA_{k\ell}(X_p a_{\ell\,j})\varphi_j\right),\,\psi\right)+\left(\Lambda^\varepsilon
\left(\sum_{\ell,j=1}^qA_{k\ell} a_{\ell\,j}(X_p \varphi_j)\right),\,\psi\right)
$$
on $U$ for $1\leq k\leq q$, which is the same as
$$
\left(\Lambda^\varepsilon\left(\sum_{\ell,j=1}^qA_{k\ell}(X_p a_{\ell\,j})\varphi_j\right),\,\psi\right)+\bigg(\Lambda^\varepsilon\left(\alpha(X_p \varphi_k)\right),\,\psi\bigg),\leqno{(5.2.1.2)}
$$
because $\sum_{j=1}^q A_{ij}a_{j\ell}=\alpha\delta_{i\ell}$ for $1\leq j,\ell\leq q$.

\medbreak We apply
$\sum_{1\leq p,k\leq q}\Gamma_{pk}$ to (5.2.1.2)
to get the estimability of
\small$$
\left(\Lambda^\varepsilon\left(\sum_{1\leq p,\ell,j,k\leq q}\Gamma_{pk}A_{k\ell}(X_p a_{\ell\,j})\varphi_j\right),\,\psi\right)+
\left(\Lambda^\varepsilon\left(\alpha\sum_{1\leq p,k\leq q}\Gamma_{pk}(X_p \varphi_k)\right),\,\psi\right)
$$
on $U$,
which is the same as
$$
\left(\Lambda^\varepsilon\left(\sum_{1\leq p,\ell,j,k\leq q}\Gamma_{pq}A_{q\ell}(X_p a_{\ell\,j})\varphi_j\right),\,\psi\right)+
\left(\sum_{1\leq p,k\leq q}\Gamma_{pq}(X_p \varphi_q),\,\Lambda^\varepsilon\left(\overline\alpha\,\psi\right)\right)
$$
up to estimable error terms from the commutators of pseudodifferential operators.

\medbreak
Since $\bar\alpha$ is a scalar-multiplier at $0$ (on account of $\alpha$ being a scalar multiplier at $0$), from the estimability of $\sum_{1\leq p,k\leq q}\Gamma_{pk}(X_p \varphi_k)$ on $U$ and the Cauchy-Schwarz inequality we conclude that
$$
\left(\sum_{1\leq p\leq r,1\leq k\leq q}\Gamma_{pk}(X_p\varphi_k),\,\Lambda^\varepsilon\left(\overline\alpha\,\psi\right)\right)
$$
is estimable on $U$.  Hence
$$
\left(\Lambda^\varepsilon\left(\sum_{1\leq p\leq r,\,1\leq\ell,j,k\leq q}\Gamma_{pk}A_{k\ell}(X_p a_{\ell\,j})\varphi_j\right),\,\psi\right)
$$
is estimable on $U$.

\medbreak We can now choose
$$
\psi=\sum_{1\leq p,\ell,j,k\leq q}\Gamma_{pk}A_{k\ell}(X_pa_{\ell\,j})\varphi_j
$$
so that
$$
\left(\Lambda^\varepsilon\left(\sum_{1\leq p\leq r,\,1\leq\ell,j,k\leq q}\Gamma_{pk}A_{k\ell}(X_p a_{\ell\,j})\varphi_j\right),\,\psi\right)
$$
is equal to
$$
\left\|\Lambda^{\frac{\varepsilon}{2}}\sum_{j=1}^q b_j\varphi_j\right\|^2
$$
This means that $\vec b$ is a vector multiplier at $0$ whose order of subellipticity is $\geq\frac{\varepsilon}{2}$.

\medbreak We now look at the two special cases.  The special case (i) follows from setting $(A_{jk})_{1\leq j,k\leq q}$ to be the adjoint matrix ${\rm adj}({\mathbf a})$ of the matrix ${\mathbf a}$ and setting $\alpha$ to be $\det({\mathbf a)}$ with $\varepsilon_2=\varepsilon_1$. The special case (ii) follows from setting $A_{jk}$ to be the Kronecker delta $\delta_{jk}$ for $1\leq j,k\leq q$ and setting $a_{jk}$ to be $\alpha\delta_{jk}$ for $1\leq j,k\leq q$ with $\varepsilon_1=\varepsilon_2$ 
Q.E.D.

\bigbreak\noindent(5.3) {\it Remark.}  Though Theorem (5.2) is presented as involving interior estimates, the same argument works also in boundary situations like the complex Neumann problem where, for the argument, integration by parts is needed only for the boundary tangential directions which do not affect the condition of the test forms to be in the domain of the actual adjoint $\bar\partial^*$ of $\bar\partial$.  The special case (ii) of Theorem (5.2), after modification for application to the situation of the complex Neumann problem for special domains, gives a procedure to generate a vector multiplier from a matrix multiplier.  In (6.5) below, computations of examples are given to show that this procedure is a new procedure of generating vector multipliers for special domains in ${\mathbb C}^n$ with $n\geq 4$.

\bigbreak\noindent(5.4) {\it Estimable Linear Combinations and Initial Multipliers.}  The goal of the techninque of multiplier ideal sheaves is to use the differential relations among the multipliers and some initial multipliers to conclude, under some geometric conditions, that the function which is identically $1$ is a scalar multiplier.  An increase in the  differential relations among the multipliers facilitates the achievement of the goal.  Theorem (5.2) uses the collection $\left(\Gamma_{k\ell}\right)_{1\leq k,\ell\leq q}$ of smooth functions on $\Omega$ to construct a new vector multiplier from a matrix multiplier.  The condition on the collection $\left(\Gamma_{k\ell}\right)_{1\leq k,\ell\leq q}$ of smooth functions on $\Omega$ is that $\sum_{1\leq k,\ell\leq q}\Gamma_{k\ell}X_k\varphi_\ell$ is estimable on $\Omega$.  For that reason we refer to the collection $\left(\Gamma_{k\ell}\right)_{1\leq k,\ell\leq q}$ of smooth functions on $\Omega$ an {\it estimable linear combination}. To facilitate the construction of new vector multipliers, we can use a family of such estimable linear combinations $\left(\Gamma_{k\ell}^{(\lambda)}\right)_{1\leq k,\ell\leq q}$ indexed by $1\leq\lambda\leq\lambda^*$ instead of a single one.  There remains the crucial question of geometric conditions to guarantee solution of the regularity problem.  This condition (which is similar to the condition of finite type for the complex Neumann problem) should be a condition on the family of estimable linear combinations $\left(\Gamma_{k\ell}^{(\lambda)}\right)_{1\leq k,\ell\leq q}$ for $1\leq\lambda\leq\lambda^*$ and the choice of initial scalar multipliers $\alpha^{(\sigma)}$ (for $1\leq\sigma\leq\sigma^*$) and initial vector multipliers ${\vec a}^{{}^{(\tau)}}$ (for $1\leq\tau\leq\tau^*$).  This question has not yet been satisfactorily answered.

\bigbreak\noindent{\bf\S 6.} {\sc New Procedure to Generate Vector Multiplier from Matrix Multiplier in Complex Neumann Problem of Special Domain}

\bigbreak We now modify the argument in the special case (ii) of Theorem (5.2) to apply to the complex Neumann problem to obtain a new procedure of generating a vector multiplier from a matrix multiplier.  This new procedure works for any bounded weakly pseudoconvex domain with smooth boundary, but we will carry out the modification only for a special domain, because the notations for a special domain have already been introduced here to make the argument for a special domain easier to present. Then we show by explicit computation for some special domains in ${\mathbb C}^4$ that this new way of generating a vector multiplier cannot be derived from Kohn's procedures in (2.1).
 
\bigbreak\noindent(6.1) {\it New Procedure to Generate Vector Multiplier from Matrix Multiplier for the Complex Neumann Problem.}  Let $\Omega$ be a special domain in ${\mathbb C}^{n+1}$ (with coordinates $w,z_1,\cdots,z_n$) defined by holomorphic functions $F_j(z_1,\cdots,z_n)$ on some open neighborhood of $\bar\Omega$, as described in (2.8.1).  For the complex Neumann problem for the special domain $\Omega$ in ${\mathbb C}^{n+1}$, the roles of the vector fields $X_1,\cdots,X_q$ are played by $\partial_j=\frac{\partial}{\partial z_j}$ for $1\leq j\leq n$ and the roles of $X_1^*,\cdots,X_q^*$ are played by $\partial_{\bar j}=\frac{\partial}{\partial\bar z_j}$ for $1\leq j\leq n$ and the role of $\Gamma_{jk}$ for $1\leq j,k\leq q$ is played by the Kronecker delta $\delta_{jk}$ for $1\leq j,k\leq n$.  Let ${\mathbf a}=\left(a_{jk}\right)_{1\leq j,k\leq n}$ be a matrix whose entry $a_{jk}$ is a holomorphic function of $z_1,\cdots,z_n$ defined on an open neighborhood of $0$ in ${\mathbb C}^n$ for $1\leq j,k\leq n$ such that each row vector $\vec a_j=(a_{jk})_{1\leq k\leq n}$ is a vector multiplier at $0$ with order of subellipticity $\geq\eta$ (for some $0<\eta\leq 1$).  Let $U$ be an open neighborhood of $0$ in ${\mathbb C}^{n+1}$ such that each $a_{jk}$, as a holomorphic function in $z_1,\cdots,z_n,w$ but independent of $w$, is defined on $U$.

\medbreak Let $\varphi=\sum_{j=1}^n\varphi_{\bar j}d\bar z_j+\hat\varphi d\bar w$ be a smooth test $(1,0)$-form on $\bar\Omega\cap U$ with compact support which is in the domain of the actual adjoint $\bar\partial^*$ of $\bar\partial$.
Let $\psi$ be a scalar function which is a linear combination of $\varphi_{\bar j}$ with smooth functions as coefficients and which we will specify more precisely later.   Let $0<\varepsilon<\frac{\eta}{2}$.  The $L^2$ inner product
$$
\left(\Lambda^{2\varepsilon}\left(\sum_{j=1}^n a_{\ell\,j}\,\varphi_{\bar j}\right),\,\bar\partial_p\psi\right)
$$
is estimable on $U$ by the Cauchy-Schwarz inequaltiy from the assumption that $\vec a_\ell=\sum_{j=1}^n a_{\ell\,j}dz_j$ is a vector-multiplier at $0$ with order of subellipticity $\geq\eta\geq 2\varepsilon$ for $1\leq\ell\leq n$ and the assumption that the $L^2$ norm of $\bar\partial_p\psi$ is estimable on $U$ for $1\leq p\leq n$ (from the estimability of $\bar\partial_p\varphi_{\bar j}$ on $U$ for any $1\leq p,\,j \leq n$).  Integration by parts applied to $\bar\partial_p$ yields the estimability of
$$
\left(\Lambda^{2\varepsilon}\left(\sum_{j=1}^n (\partial_p a_{\ell\,j})\varphi_{\bar j}\right),\,\psi\right)
+\left(\Lambda^{2\varepsilon}\left(\sum_{j=1}^n a_{\ell\,j}(\partial_p \varphi_{\bar j})\right),\,\psi\right)\leqno{(6.1.1)}
$$
on $U$ after we take care of the error terms from the commutators of operators in the standard way.

\bigbreak Let $\left(A_{q\ell}\right)_{1\leq q,\ell\leq n}$ be the adjoint matrix  of ${\mathbf a}$ so that $\sum_{j=1}^n A_{ij}a_{jk}=(\det{\mathbf a})\delta_{ik}$ for $1\leq i,k\leq n$ (where $\delta_{jk}$ is the Kronecker delta).  Now for $1\leq p\leq n$ we apply $\sum_{\ell=1}^n A_{p\ell}$ to (6.1.1) to get the estimability of
$$
\left(\Lambda^{2\varepsilon}\left(\sum_{\ell,j=1}^nA_{p\ell}(\partial_p a_{\ell\,j})\varphi_{\bar j}\right),\,\psi\right)+\left(\Lambda^{2\varepsilon}
\left(\sum_{\ell,j=1}^n A_{p\ell} a_{\ell\,j}(\partial_p \varphi_{\bar j})\right),\,\psi\right)
$$
on $U$, which is the same as
$$
\left(\Lambda^{2\varepsilon}\left(\sum_{\ell,j=1}^n A_{p\ell}(\partial_p a_{\ell\,j})\varphi_{\bar j}\right),\,\psi\right)+\bigg(\Lambda^{2\varepsilon}\left((\det{\mathbf a})(\partial_p\varphi_{\bar p})\right),\,\psi\bigg),\leqno{(6.1.2)}
$$
because $\sum_{j=1}^n A_{ij}a_{jk}=(\det{\mathbf a})\delta_{ik}$ for $1\leq i,k\leq n$.

\bigbreak We now sum (6.1.2) over $1\leq p\leq n$ to get the estimability of
$$
\left(\Lambda^{2\varepsilon}\left(\sum_{p,\ell,j=1}^n A_{p\ell}(\partial_p a_{\ell\,j})\varphi_{\bar j}\right),\,\psi\right)+\left(\Lambda^{2\varepsilon}\left((\det{\mathbf a})\left(\sum_{p=1}^n\partial_p\varphi_{\bar p}\right)\right),\,\psi\right)
$$
on $U$, which is the same as
$$
\left(\Lambda^{2\varepsilon}\left(\sum_{p,\ell,j=1}^nA_{p\ell}(\partial_p a_{\ell\,j})\varphi_{\bar j}\right),\,\psi\right)+\left(\sum_{p=1}^n\partial_p\varphi_{\bar p},\,\Lambda^{2\varepsilon}\left(\overline{(\det{\mathbf a})}\psi\right)\right).
$$
As a determinant whose rows are vector-multipliers with order of subellipticity $\geq\eta\geq 2\varepsilon$, the determinant $\det({\mathbf a})$ (as well as its complex-conjugate) is a scalar-multiplier with order of subellipticity $\geq\eta\geq 2\varepsilon$.  The function
$\Lambda^{2\varepsilon}\left(\overline{(\det{\mathbf a})}\psi\right)$ is estimable on $U$, because $\psi$ is a linear combination of $\varphi_{\bar j}$ for $1\leq j\leq n$ with smooth functions as coefficients.
From the estimability of $\sum_{p=1}^n\partial_p\varphi_{\bar p}$ on $U$ and the Cauchy-Schwarz inequality we conclude that
$$
\left(\Lambda^{2\varepsilon}\left(\sum_{p,\ell,j=1}^n A_{p\ell}(\partial_p a_{\ell\,j})\varphi_{\bar j}\right),\,\psi\right)
$$
is estimable on $U$. We now choose
$$
\psi=\sum_{p,\ell,j=1}^nA_{p\ell}(\partial_p a_{\ell\,j})\varphi_{\bar j}
$$
so that
$$
\left(\Lambda^{2\varepsilon}\left(\sum_{p,\ell,j=1}^nA_{p\ell}(\partial_p a_{\ell\,j})\varphi_{\bar j}\right),\,\psi\right)_{L^2(\Omega)}
$$
is equal to $\left\|\Lambda^\varepsilon\psi\right\|^2_{L^2(\Omega)}$.
This means that the $(1,0)$-form
$$
\sum_{j=1}^n\left(\sum_{1\leq p,\ell\leq n}A_{p\ell}(\partial_p a_{\ell\,j})\right)dz_j
$$
is a vector-multiplier with order of subellitpicity $\geq\frac{\eta}{2}$..  This is a new process of producing vector-multipliers from a matrix of vector-multipliers.  We now summarize the result in the following theorem.

\bigbreak\noindent(6.2) {\it Theorem.}  Let $\Omega$ be a special domain in ${\mathbb C}^{n+1}$ (with coordinates $w,z_1,\cdots,z_n$) defined by (2.8.1.1).  Assume that $0$ is a boundary point of $\Omega$.  Let ${\mathbf a}=\left(a_{jk}\right)_{1\leq j,k\leq n}$ be a matrix whose entry $a_{jk}$ is a holomorphic function of $z_1,\cdots,z_n$ defined on an open neighborhood of $0$ in ${\mathbb C}^n$ for $1\leq j,k\leq n$ such that each row vector $\vec a_j=(a_{jk})_{1\leq k\leq n}$ is a vector multiplier at $0$ with order of subellipticity $\geq\eta$ (for some $0<\eta\leq 1$).  Let ${\rm adj}({\mathbf a})$ be the adjoint matrix of ${\mathbf a}$.  Then the holomorphic $(1,0)$-form
$$
\sum_{j=1}^n\left(\sum_{1\leq p,\ell\leq n}\left({\rm adj}({\mathbf a})\right)_{p\ell}(\partial_p a_{\ell\,j})\right)dz_j\leqno{(6.2.1)}
$$
is a vector multiplier at $0$ with order of subellipticity $\geq\frac{\eta}{2}$, where $\left({\rm adj}({\mathbf a})\right)_{p\ell}$ is the entry of ${\rm adj}({\mathbf a})$ in the $p$-th row and the $\ell$-th column.

\bigbreak\noindent(6.3) {\it Comparison with Known Procedure of Generating Vector Multiplier from Matrix Multiplier.}   In the case of a special domain, the known procedure (B)(i)(ii) in (2.1) to generate a vector multiplier from a given matrix multiplier ${\mathbf a}$ is to first use (B)(ii) in (2.1) to get the determinant $\det({\mathbf a})$ of ${\mathbf a}$ as a scalar multiplier and then use (B)(i) in (2.1) to get the $(1,0)$-form $\partial(\det({\mathbf a}))$ as a vector multiplier.  Here we use the same notations as in (6.2).  

\medbreak We would like to compare $\partial(\det({\mathbf a}))$ with the vector multiplier from Theorem (6.2).  Since 
$\sum_{j=1}^n\left({\rm adj}({\mathbf a})\right)_{kj}a_{jk}=\det({\mathbf a})$ for any $1\leq k\leq n$, we have
$$
\partial(\det{\mathbf a})=\sum_{j=1}^n(\partial\left({\rm adj}({\mathbf a})\right)_{kj})a_{jk}+\sum_{j=1}^n\left({\rm adj}({\mathbf a})\right)_{kj}(\partial a_{jk})\leqno{(6.3.1)}
$$
which is different from the vector multiplier (6.2.1)
where the index $j$ of $a_{\ell\,j}$ is used as the index for the component of the vector multiplier instead of the subscript $j$ of $\partial_j$ which is the index for the component
of the vector multiplier (6.3.1).

\medbreak As shown in the computations given below in (6.4) and (6.5), it turns out that in the case of $n=2$ the old procedure (B)(i)(ii) in (2.1) produces the same result as the new procedure given in Theorem (6.2), but in the case of $n\geq 3$ the new procedure indeed gives some new vector multipliers different from those produced by the procedures B(i)(ii).

\bigbreak\noindent(6.4) {\it New Procedure Gives No New Vector Multipliers for $3$-Dimensional Special Domain.}  We explicitly compute (6.2.1) and (6.3.1) in the case of a special domain $\Omega$ in ${\mathbb C}^3$ to determine whether the result (6.2.1) from the new procedure is different from the result (6.3.1) from the old procedure.  We need only consider holomorphic functions and holomorphic $1$-forms on ${\mathbb C}^2$ as scalar and vector multipliers.  Let $a_{j1}dz_1+a_{j2}dz_2$ as two holomorphic $1$-forms which are vector multipliers.  For the matrix multiplier ${\mathbf a}=\left(a_{jk}\right)_{1\leq j,k\leq 2}$, the adjoint matrix ${\rm adj}({\mathbf a})$ is
$$
\left(A_{jk}\right)_{1\leq j,k\leq 2}=\left(\begin{matrix}a_{22}&-a_{12}\cr -a_{21}& a_{11}\end{matrix}\right).
$$
The vector multiplier generated by the new procedure is $b_1dz_1+b_2dz_2$ with
$$
\begin{aligned}b_j&=A_{11}\partial_1 a_{1j}+A_{12}\partial_1 a_{2j}+A_{21}\partial_2 a_{1j}+A_{22}\partial_2 a_{2j}\cr
&=a_{22}\partial_1 a_{1j}-a_{12}\partial_1 a_{2j}-a_{21}\partial_2 a_{1j}+a_{11}\partial_2 a_{2j}.\cr
\end{aligned}
$$
On the other hand,
$$
\partial_j(\det{\mathbf a})=a_{11}\partial_j a_{22}+a_{22}\partial_j a_{11}-a_{12}\partial_j a_{21}-a_{21}\partial_j a_{12}
$$
so that the difference of their dot products with a test $(0,1)$-form $\varphi$ is
$$
\begin{aligned}&(b_1dz_1+b_2dz_2)\cdot\varphi-(\partial(\det{\mathbf a}))\cdot\varphi\cr
&=\Big(a_{22}(\partial_1 a_{11})\varphi_{\bar 1}-a_{12}(\partial_1 a_{21})\varphi_{\bar 1}-a_{21}(\partial_2 a_{11})\varphi_{\bar 1}+a_{11}(\partial_2 a_{21})\varphi_{\bar 1}\cr
&+a_{22}(\partial_1 a_{12})\varphi_{\bar 2}-a_{12}(\partial_1 a_{22})\varphi_{\bar 2}-a_{21}(\partial_2 a_{12})\varphi_{\bar 2}+a_{11}(\partial_2 a_{22})\varphi_{\bar 2}\Big)\cr
&-\Big(a_{11}(\partial_1a_{22})\varphi_{\bar 1}+a_{22}(\partial_1 a_{11})\varphi_{\bar 1}-a_{12}(\partial_1 a_{21})\varphi_{\bar 1}-a_{21}(\partial_1 a_{12})\varphi_{\bar 1}\cr
&+a_{11}(\partial_2a_{22})\varphi_{\bar 2}+a_{22}(\partial_2 a_{11})\varphi_{\bar 2}-a_{12}(\partial_2 a_{21})\varphi_{\bar 2}-a_{21}(\partial_2 a_{12})\varphi_{\bar 2}\Big)\cr
&=(\partial_2a_{21}-\partial_1a_{12})(a_{11}\varphi_{\bar 1}+a_{12}\varphi_{\bar 2})+(\partial_1a_{12}-\partial_2a_{21})(a_{21}\varphi_{\bar 1}+a_{22})\varphi_{\bar 2}\cr
\end{aligned}
$$
with the cancellation of the four pairs of terms ((1,1),(3,2)), ((1,2),(3,3)), ((2,3),(4,4)), ((2,4),(4,1)) where $(j,k)$ refers to the term on the $j$-th row in the $k$-th position in the array of terms in four rows of four terms each.  The difference is a linear combination of $a_{j1}\varphi_{\bar 1}+a_{j2}\varphi_{\bar 2}$ for $j=1,2$ with smooth coefficients and is known to be estimable, because $a_{j1}dz_1+a_{j2}dz_2$ is a vector multiplier for $j=1,2$.  It means that in the case of a special domain in ${\mathbb C}^3$ the new procedure does not give any new multipliers.

\bigbreak\noindent(6.5) {\it New Procedure Gives More Vector Multipliers for $4$-Dimensional Special Domain.}   The new procedure of generating vector multipliers already gives vector multiplies different from those generated by the procedure B(i)(ii) in (2.1) in the case of special domain in ${\mathbb C}^4$, as shown in the following computation.  Consider the upper triangular matrix
$$
{\mathbf a}=\left(\begin{matrix}a_{11}&\xi&0\cr
0&a_{22}&\eta\cr
0&0&a_{33}\cr\end{matrix}\right)
$$
whose three row vectors are vector multipliers which are holomorphic in the variables $z_1,z_2,z_3$.  Its adjoint matrix ${\rm adj}({\mathbf a})$, as the inverse of ${\mathbf a}$ times its determinant, is
$$
\left(A_{jk}\right)_{1\leq j,k\leq 3}=\left(\begin{matrix}a_{22}a_{33}&-a_{33}\xi&\xi\eta\cr
0&a_{11}a_{33}&-a_{11}\eta\cr
0&0&a_{11}a_{22}\cr\end{matrix}\right).
$$
We now compare the vector multiplier $\partial(\det\,{\mathbf a})$ with the vector multiplier $\sum_{j=1}^3 b_jdz_j$ generated by the new procedure by taking the difference of their dot products with a test function $\varphi=\sum_{j=1}^3\varphi_{\bar j}d\bar z_j+\hat\varphi dw$ to get
$$\begin{aligned}&\left(\sum_{j=1}^3 b_jdz_j\right)\cdot\varphi-\partial(\det\,{\mathbf a})\cdot\varphi\cr
&=\sum_{1\leq j<k\leq 3}A_{jk}(\partial_\nu a_{jk})\varphi_{\bar\nu}\cr
&=\sum_{\nu=1}^3\left(-a_{33}\xi\partial_\nu\xi-a_{11}\eta\partial_\nu\eta\right)\varphi_{\bar\nu}\cr
&=-\frac{1}{2}\sum_{\nu=1}^3\left(a_{33}\partial_\nu(\xi^2)-a_{11}\partial_\nu(\eta^2)\right)\varphi_{\bar\nu}.\cr
\end{aligned}
$$
The assumption that ${\mathbf a}$ is a matrix multiplier gives us only the estimability of
$$
a_{11}\varphi_{\bar 1}+\xi\varphi_{\bar 2},\quad a_{22}\varphi_{\bar 2}+\eta\varphi_{\bar 3},\quad a_{33}\varphi_{\bar 3},
$$
from which we cannot conclude the estimability of
$$
-\frac{1}{2}\sum_{\nu=1}^3\left(a_{33}\partial_\nu(\xi^2)-a_{11}\partial_\nu(\eta^2)\right)\varphi_{\bar\nu}
$$
for arbitrary functions $\xi$ and $\eta$, because even in the special case of $\varphi_{\bar 1}=\varphi_{\bar 2}=0$, the estimability  of
$$
-\frac{1}{2}\sum_{\nu=1}^3\left(a_{33}\partial_\nu(\xi^2)-a_{11}\partial_\nu(\eta^2)\right)\varphi_{\bar\nu}
$$
would mean the estimability of $a_{11}\partial_3(\eta^2)\varphi_{\bar 3}$ for arbitrary $a_{11}$ and $\eta$, which cannot be derived from the estimability of $a_{33}\varphi_{\bar 3}$.  This shows that the new procedure gives more vector multipliers.

\bigbreak\noindent{\sc References}

\medbreak\noindent[Catlin-D'Angelo2010] David Catlin and John P. D'Angelo,
Subelliptic estimates. In: {\it Several Complex Variables and Connections with PDE Theory and Geometry (Proc. 2008 Fridbourg Conference for Linda Rothschild)}. ed. Ebenfelt {\it et al}. Trends in Mathematics. Birkh\"auser/Springer Basel AG, Basel, 2010, pp.75 -- 94.

\medbreak\noindent[D'Angelo1979] John P. D'Angelo, Finite type
conditions for real hypersurfaces. {\it J. Differential Geom.} \textbf{14} (1979),
59 -- 66.

\medbreak\noindent[D'Angelo1982] John P. D'Angelo,
Real hypersurfaces, orders of contact, and applications.
{\it Ann. of Math.} \textbf{115} (1982),  615 –- 637.

\medbreak\noindent[Diederich-Fornaess1978] Klas Diederich and John-Erik Fornaess,
Pseudoconvex domains with real-analytic boundary. {\it Ann. of
Math.} \textbf{107} (1978), 371 -- 384.

\medbreak\noindent[Donaldson1985] Simon Donaldson,
Anti self-dual Yang-Mills connections over complex algebraic surfaces and stable vector bundles.
{\it Proc. London Math. Soc.} \textbf{50} (1985), 1 –- 26.

\medbreak\noindent[H\"ormander1967] Lars H\"ormander,
Hypoelliptic second order differential equations.
{\it Acta Math.} \textbf{119} (1967), 147 –- 171.

\medbreak\noindent[Kohn1972]
Joseph J. Kohn,
Boundary behavior of $\bar\partial$ on weakly pseudo-convex manifolds of dimension two.
Collection of articles dedicated to S. S. Chern and D. C. Spencer on their sixtieth birthdays.
{\it J. Differential Geometry} \textbf{6} (1972), 523 –- 542.

\medbreak\noindent[Kohn1977]
Joseph J. Kohn,
Sufficient conditions for subellipticity on weakly pseudo-convex domains.
{\it Proc. Nat. Acad. Sci. U.S.A.} \textbf{74} (1977),  2214 –- 2216.

\medbreak\noindent[Kohn1979] Joseph J. Kohn, Subellipticity of the
$\bar \partial $-Neumann problem on pseudo-convex domains:
sufficient conditions. {\it Acta Math.} \textbf{142} (1979), 79 -- 122.

\medbreak\noindent[Kohn-Nirenberg1965]
Joseph J. Kohn and Louis Nirenberg,
Non-coercive boundary value problems.
{\it Comm. Pure Appl. Math.} \textbf{18} (1965), 443 –- 492.

\medbreak\noindent[Nadel1990] Alan Michael Nadel,
Multiplier ideal sheaves and Kähler-Einstein metrics of positive scalar curvature.
{\it Ann. of Math.} \textbf{132} (1990), 549 –- 596.

\medbreak\noindent[Siu2010] Yum-Tong Siu,
Effective Termination of Kohn's Algorithm for
Subelliptic Multipliers, {\it Pure and Applied Mathematics Quarterly} \textbf{6} (2010)
(Special Issue in honor of Joseph J. Kohn), 1169 -- 1241. (arXiv:0706.4113)

\medbreak\noindent[Skoda1972] Henri Skoda,
Application des techniques $L^2$ \`a la th\'eorie des id\'eaux d'une alg\`ebre de fonctions holomorphes avec poids.
{\it Ann. Sci. \'Ecole Norm. Sup.} \textbf{5} (1972), 545 -- 579.

\bigbreak\noindent{\it Author's address:} Department of Mathematics, Harvard University, Cambridge, MA 02238, U.S.A.
\end{document}